\documentclass[final]{siamart}
\usepackage{lipsum}
\usepackage{amsfonts}
\usepackage{graphicx}
\graphicspath{{./images/}}
\usepackage{epstopdf}
\usepackage{algorithmic}
\ifpdf
  \DeclareGraphicsExtensions{.eps,.pdf,.png,.jpg}
\else
  \DeclareGraphicsExtensions{.eps}
\fi

\usepackage{labtex}

\newcommand{\TheTitle}{Linearly Solvable Stochastic Control Lyapunov Functions} 
\newcommand{\TheAuthors}{Y.~P. Leong, M.~B. Horowitz, and J.~W. Burdick}

\headers{\TheTitle}{\TheAuthors}

\title{{\TheTitle}\thanks{A short version of this work appeared in \cite{ypl2015linearsoc}.}}

\author{
  Yoke Peng Leong\thanks{Control and Dynamical Systems, California Institute of Technology, 
Pasadena, CA
    (\email{ypleong@caltech.edu}, \email{mhorowit@caltech.edu}).}
  \and
  Matanya B. Horowitz\footnotemark[2]
  \and
  Joel W. Burdick \thanks{Mechanical Engineering, California Institute of Technology, Pasadena, CA (\email{jwb@robotics.caltech.edu}).}
}

\usepackage{amsopn}

\newsiamremark{rem}{Remark}
\newsiamthm{ass}{Assumption}

\usepackage[tight,footnotesize]{subfigure}

\begin{document}
\maketitle

\begin{abstract}                          % Abstract of not more than 250 words.
This paper presents a new method for synthesizing stochastic control Lyapunov functions for a class of nonlinear stochastic control systems. The technique relies on a transformation of the classical nonlinear Hamilton-Jacobi-Bellman partial differential equation to a linear partial differential equation for a class of problems with a particular constraint on the stochastic forcing. This linear partial differential equation can then be relaxed to a linear differential inclusion, allowing for relaxed solutions to be generated using sum of squares programming. The resulting relaxed solutions are in fact viscosity super/subsolutions, and by the maximum principle are pointwise upper and lower bounds to the underlying value function, even for coarse polynomial approximations. Furthermore, the pointwise upper bound is shown to be a stochastic control Lyapunov function, yielding a method for generating nonlinear controllers with pointwise bounded distance from the optimal cost when using the optimal controller. These approximate solutions may be computed with non-increasing error via a hierarchy of semidefinite optimization problems. Finally, this paper develops a-priori bounds on \emph{trajectory suboptimality} when using these approximate value functions, as well as demonstrates that these methods, and bounds, can be applied to a more general class of nonlinear systems not obeying the constraint on stochastic forcing. Simulated examples illustrate the methodology.
\end{abstract}



%%% ------------------------------------------------------------------------------------------- %%%
\section{Introduction} \label{sec:intro}

The study of system stability is a central theme of control engineering. A primary tool for such studies is Lyapunov theory, wherein an energy-like function is used to show that some measure of distance from a stability point decays over time.  The construction of Lyapunov functions that certify system stability advanced considerably with the introduction of Sums of Squares (SOS) programming, which has allowed for Lyapunov functions to be synthesized for both polynomial systems
\cite{Parrilo:2003fh} and more general vector fields \cite{papachristodoulou2005analysis}.

To address the more challenging problem of stabilization, rather than the analysis of an existing closed loop system, it is possible to generalize Lyapunov functions to incorporate control inputs. The existence of a Control Lyapunov Function (CLF) (see \cite{krstic1995nonlinear,freeman1996control,sontag1983lyapunov}) is sufficient for the construction of a stabilizing controller. However, the synthesis of a CLF for general systems remains an open question. Unfortunately, the SOS-based methods cannot be naively extended to the generation of CLF solutions, due to the bilinearity between the Lyapunov function and control input.

% However, for several large and important classes of systems, CLFs are in fact known and may be used for stabilization, with a review of the theory available in \cite{sontag1983lyapunov}. The drawback is that these CLFs are hand-constructed and may be shown to be arbitrarily suboptimal. 

Due to the lack of a general CLF synthesis technique, an alternative is the use of Receding Horizon Control (RHC), which allows for the incorporation of optimality criteria. Euler-Lagrange equations are used to construct a locally optimum trajectory \cite{Primbs:1999tt}, and stabilization is guaranteed by constraining the terminal cost in the RHC problem to be a CLF. Suboptimal CLFs have found extensive use, with applications in legged locomotion \cite{kolathaya2014exponential} and distributed control \cite{ogren2001control}. Adding stochasticity to the governing dynamics compounds the difficulties of constructing Lyapunov functions \cite{deng1997stochastic,florchinger1997feedback}.
%\ypl{We should probably focus on the stochastic literatures instead of the deterministic ones.}
A complementary area in control engineering is the study of the Hamilton-Jacobi-Bellman (HJB) equation that governs the optimal control of a system. Methods to calculate the solution to the HJB equation via semidefinite programming have been proposed previously by Lasserre \emph{et al.} \cite{lasserre2008nonlinear}.
 % In their work, the solution and the optimality conditions are integrated against monomial test functions, producing an infinite set of moment constraints. By truncating to any finite list of monomials, the optimal control problem is reduced to a semidefinite optimization problem. 
 The method is quite general, applicable to any system with polynomial nonlinearities. 
% These moment techniques are also intimately related to sum of squares programming, and it can be shown that the two problems are in fact convex duals of one another \cite{lasserre2008nonlinear}. 

In this work, we propose an alternative line of study based on the linear structure of a particular form of the HJB equation. Since the late 1970s, Fleming \cite{fleming1982logarithmic}, Holland \cite{holland1977new} and other researchers thereafter \cite{dai1991stochastic,filliger2005relative} have made connections between stochastic optimal control and reaction-diffusion equation through a logarithmic transformation. Recently, when studying stochastic control using the HJB equation, Kappen \cite{Kappen:2005kb} and Todorov \cite{Todorov:2009wja} discovered that particular assumptions on the structure of a dynamical system, given the name \emph{linearly solvable} systems, allows a logarithmic transformation of the optimal control equation to a linear partial differential equation form. The linearity of this class of problems has given rise to a growing body of research, with an overview available in \cite{Dvijotham:2012tv}. Kappen's work focused on calculating solutions via \emph{path integral} techniques. Todorov began with the analysis of particular Markov decision processes, and showed the connection between the two paradigms. This work was built upon by Theodorou et al. \cite{Theodorou:2011uz} into the Path Integral framework in use with Dynamic Motion Primitives. These results have been developed in many compelling directions \cite{Stulp:2012jb,Dvijotham:2012tv,wiegerinck2012stochastic,rutquist2014hjb}.

This paper combines these previously disparate fields of linearly solvable optimal control and Lyapunov theory, and provides a systematic way to construct stabilizing controllers with guaranteed performance. The result is a hierarchy of SOS programs that generates stochastic CLFs (SCLF) for arbitrary \emph{linearly solvable} systems. Such an approach has many benefits. First and foremost, this approach generates stabilizing controllers for an important class of nonlinear, stochastic systems even when the optimal controller is not found. We prove that the approximate solutions generated by the SOS programs are pointwise upper and lower bounds to the true solutions. In fact, the upper bound solutions are SCLFs which can be used to construct stabilizing controllers, and they bound the performance of the system when they are used to construct suboptimal controllers. Existing methods for the generation of SCLFs do not have such performance guarantees. 
% Hence, when users are limited by their computing resources, they could use the best suboptimal solutions they can generate to stabilize the system with guaranteed performance. 
Additionally, we demonstrate that, although the technique is based on linear solvability, it may be readily extended to more general systems, including deterministic systems, while inheriting the same performance guarantees. 

% Furthermore, the HJB solution is a globally optimal controller, incorporating all potential initial system states and the system dynamics. 

% Here, we propose polynomial candidate approximate solutions to the HJB, extending recently developed tools in polynomial optimization to a new class of problems. It is already known that the solution to the deterministic HJB is in fact a CLF \cite{Primbs:1999tt}, and similarly the stochastic HJB is a SCLF. This paper shows that our approximated value function solutions are also SCLFs. 

A preliminary version of this work appeared in \cite{horowitz2014semidefinite} and \cite{horowitz2014admm}, where the use of sum of squares programming for solving the HJB were first considered. This paper builds on this recent body of research, studying the stabilization and optimality properties of the resulting solutions. These previous works focused on path planning, rather than stabilization, and did not include the stability analysis or suboptimality guarantees presented in this paper. A short version of this work appeared in \cite{ypl2015linearsoc} which included less details and did not include the extension in \cref{sec:extensions}.

The rest of this paper is organized as follows. \Cref{sec:background} reviews linearly solvable HJB equations, SCLFs, and SOS programming. \Cref{sec:problem} introduces a relaxed formulation of the HJB solutions which is efficiently computable using the SOS methodology. \Cref{sec:analysis} analyzes the properties of the relaxed solutions, such as approximation errors relative to the exact solutions. This section shows that the relaxed solutions are SCLFs, and that the resulting controller is stabilizing. The upper bound solution is also shown to bound the performance when using the suboptimal controller. \Cref{sec:extensions} summarizes an extension of the method to approximate optimal control problems which are not linearly solvable. Two examples are presented in \cref{sec:simulation} to illustrate the optimization technique and its performance. \Cref{sec:conclusion} summarizes the findings of this work and discusses future research directions.

%%% ------------------------------------------------------------------------------------------- %%%
\section{Backgrounds} \label{sec:background} 

This section briefly describes the paper's notation and reviews necessary background on the linear HJB equation, SCLFs, and SOS programming.

%%% ------------------------------------------------------------------------------------------- %%%
\subsection{Notation} \label{sec:notation}

\Cref{tab:notation} summarizes the notation of different sets appearing in the paper.
\begin{table}[!ht] 
	\centering
	\caption{Set notation}
	\label{tab:notation}
\begin{tabular}{p{0.11\textwidth}p{0.8\textwidth}}
	\hline
	Notation & Definition \\
	\hline 
        $\mathbb{Z}_+$ & All positive integers \\
	$\mathbb{R}$ & All real numbers \\
	$\mathbb{R}_+$ & All nonnegative real numbers \\
	$\mathbb{R}^n$ & All $n$-dimensional real vectors \\
	$\mathbb{R}[x]$ & All real polynomial functions in $x$ \\
	$\mathbb{R}^{n \times m}$ & All $n \times m$ real matrices\\
	$\mathbb{R}^{n \times m}[x]$ & All $M \in \mathbb{R}^{n \times m}$ such that $M_{i,j} \in \mathbb{R}[x]~ \forall~i,j$\\
	$\mathcal{K}$ & All continuous nondecreasing functions $\mu: \mathbb{R}_+ \to \mathbb{R}_+$ such that  $\mu(0) = 0$, $\mu(r ) > 0$ if $r > 0$, and $\mu(r) \ge \mu(r')$ if $r > r'$ \\
	$\mathcal{C}^{k,k'}$ & All functions $f$ such that $f$ is $k$-differentiable with respect to the first argument and $k'$-differentiable with respect to the second argument \\
	\hline
\end{tabular}
\end{table}

A compact domain in $\mathbb{R}^n$ is denoted as $\Omega$ where $\Omega \subset \mathbb{R}^n$, and its boundary is denoted as $\partial \Omega$. A domain $\Omega$ is a \emph{basic closed semialgebraic} set if there exists $g_i(x) \in \mathbb{R}[x]$ for $ i = 1,2,\ldots,m$ such that $\Omega = \{x \mid g_i(x) \ge 0 ~ \forall i = 1,2,\ldots,m\}$.  %The set of integer is denoted as $\mathbb{Z}$, and the set of positive integer is denoted as $\mathbb{Z}_+$.

A point on a trajectory, $x(t)\in\mathbb{R}^n$, at time $t$ is denoted $x_{t}$, while the segment of this trajectory over the interval $[t,T]$ is denoted by $x_{t:T}$.  %The full trajectory from time 0 to a final time $T$ is denoted as $\mathbf{x}$.

Given a polynomial $p(x)$, $p(x)$ is positive on domain $\Omega$ if $p(x) >0 ~\forall x \in \Omega$, $p(x)$ is nonnegative on domain $\Omega$ if $p(x) \ge 0 ~\forall x \in \Omega$, and $p(x)$ is positive definite on domain $\Omega$ where $0 \in \Omega$, if $p(0) = 0$ and $p(x) >0$ for all $x \in \Omega \backslash \{0\}$.

If it exists, the infinity norm of a function is defined as $\norm{f}_\infty = \sup_x |f(x)|$ for $x\in\Omega$.
To improve readability, a function, $f(x_1,\ldots,x_n)$, is abbreviated as $f$ when the arguments of the function are clear from the context.

%%% ------------------------------------------------------------------------------------------- %%%
\subsection{Linear Hamilton-Jacobi-Bellman (HJB) Equation}\label{sec:hjb}

Consider the following affine nonlinear dynamical system,
  \begin{equation} \label{eq:stochastic-dynamics}
    dx_t = \left(f(x_t) + G(x_t)u_t\right)dt + B(x_t) \, d\omega_t
  \end{equation}
where $x_t\in\Omega$ is the state at time $t$ in a compact domain $\Omega \subset \mathbb{R}^{n}$, and $u_t \in \mathbb{R}^{m}$ is the control input, $f(x)\in \mathbb{R}^n[x]$, $G(x)\in \mathbb{R}^{n \times m}[x]$, and $B(x)\in \mathbb{R}^{n \times l}[x] $ are real polynomial functions of the state variables $x$, and $\omega_{t} \in \mathbb{R}^{l}$ is a vector consisting of Brownian motions with covariance $\Sigma_\epsilon$, i.e., $\omega_t$ has independent increments with $\omega_{t}-\omega_{s}\sim\mathcal{N}(0,\Sigma_\epsilon(t-s))$, for $\mathcal{N}\left(\mu,\sigma^{2}\right)$, a normal distribution. The domain $\Omega$ is assumed to be a basic closed semialgebraic set defined as $\Omega = \{x \mid  g_i(x) \in \mathbb{R}[x], g_i(x) \ge 0 ~ \forall i = 1,2,\ldots,m\}$. %Extensions to non-polynomial functions are discussed in \cref{sec:npoly}. 
Without loss of generality, let $0 \in \Omega$ and $x=0$ be the equilibrium point, whereby $f(0) = 0$, $G(0) = 0$ and $B(0) = 0$. 

%\ypl{What does $\mathbb{E}_{\omega_t}$ mean?}
The goal is to minimize the following functional,
  \begin{gather} \label{eq:cost-functional}
       J(x,u)=\mathbb{E}_{\omega_t}\left[\phi(x_{T})+\int_{0}^{T} q(x_t)+\frac{1}{2}u_{t}^{T}Ru_{t} ~dt\right]
  \end{gather}
subject to \cref{eq:stochastic-dynamics}, where $\phi \in \mathbb{R}[x]$, $\phi: \Omega \to \mathbb{R}_+$ represents a state-dependent terminal cost, $q \in \mathbb{R}[x]$, $q: \Omega \to \mathbb{R}_+ $ is state dependent cost, and $R \in \mathbb{R}^{m \times m}$ is a positive definite matrix. The final time, $T$, unknown a priori, is the time at which the system reaches the domain boundary or the origin. This problem is generally called the \emph{first exit} problem. The expectation $\mathbb{E}_{\omega_t}$ is taken over all realizations of the noise $\omega_t$. For stability of the resultant controller to the origin, $q(\cdot)$ and $\phi(\cdot)$ are also required to be positive definite functions. %Note that the control effort enters the equation quadratically. 

The solution to this minimization problem is known as the \emph{value function}, $V: \Omega \to \mathbb{R}_+$, where beginning from an initial point
$x_t$ at time $t$
  \begin{equation} \label{eq:value-def}
   V\left(x_t\right)=\min_{u_{t:T}} J\left(x_{t:T},u_{t:T}\right).
  \end{equation}

Based on dynamic programming arguments \cite[Ch.~III.7]{Fleming:2006tl}, the associated HJB equation is a nonlinear, second order partial differential equation (PDE):
  % \begin{equation} \label{eq:original-hjb}
  %      0 = \min_{u}\left(q+\frac{1}{2}u^{T}Ru+\left(\nabla_{x}V\right)^{T}(f + G u)
  %       +\frac{1}{2}Tr\left(\left(\nabla_{xx}V\right)B\Sigma_{\epsilon}B^{T}\right)\right)
  % \end{equation}
  \begin{equation}  \label{eq:hjb-pde-value}
     0 =  q+\left(\nabla_{x}V\right)^{T}f-\frac{1}{2}
           \left(\nabla_{x}V\right)^{T}GR^{-1}G^{T}\left(\nabla_{x}V\right)
                    +  \frac{1}{2}Tr\left(\left(\nabla_{xx}V\right)B\Sigma_{\epsilon}B^{T}\right)
  \end{equation}
  with boundary condition $V(x) = \phi(x)$. For the stabilization problem on a compact domain, it is appropriate to set the boundary condition to be $\phi(x)=0$ for $x=0$, indicating zero cost accrued for achieving the origin, and $\phi(x)>0$ for $x \in \partial \Omega \setminus \{0\}$. In practice, $\phi(x)$ at the exterior boundary is usually chosen to be a large number that depends on the given application to impose large penalty for exiting the predefined domain. The optimal control effort, $u^*$, is given by
  \begin{equation}\label{eq:optimal-u}
        u^{*}=-R^{-1}G^{T}\nabla_x V.
  \end{equation}
% Substituting the optimal control, $u^{*}$, into \cref{eq:original-hjb} yields the following
% nonlinear, second order partial differential equation (PDE):

%State constraints enforced as infinite penalties, made possible via the logaritheoremic transform shown shortly, are also discussed in \cite{wik2010state}. %\ypl{What would be the right state BC? Normally, HJB has no state boundary... I think. If we set the boundary to be infinity, then Theorem \cref{thm:cost-upper} does not has a meaningful upper bound as it will be infinity.}

In general, solving \eqref{eq:hjb-pde-value} is difficult due to its nonlinearity. But, with the assumption that there
exists a $\lambda > 0$, a control penalty cost $R$ in \eqref{eq:cost-functional} satisfying
the equation
  \begin{equation} \label{eq:noise-assumption}
    \lambda G(x_t)R^{-1}G(x_t)^{T}=B(x_t)\Sigma_{\epsilon}B(x_t)^{T}\triangleq \Sigma(x_t) \triangleq  \Sigma_{t}
   \end{equation}
and using the logarithmic transformation 
  \begin{equation} \label{eq:log-transform}
      V=-\lambda\log\Psi,
  \end{equation}
it is possible \cite{EvangelosIFAC11,Todorov:2009wja,Kappen:2005bn}, after substitution and simplification, to obtain the following {\em linear}
PDE from \cref{eq:hjb-pde-value}
  \begin{gather} 
   0 = -\frac{1}{\lambda}q\Psi+f^{T}(\nabla_{x}\Psi)
      + \frac{1}{2}Tr\left(\left(\nabla_{xx}\Psi\right)\Sigma_{t}\right) \quad x \in \Omega \label{eq:hjb-pde}
      \\
   \Psi(x) = e^{- \frac{\phi(x)}{\lambda}}  \quad x \in \partial \Omega.\nonumber 
  \end{gather}

This transformation of the value function has been deemed the \emph{desirability} function \cite[Table 1]{Todorov:2009wja}. 
% Similar arguments can be used to construct desirability functions for the other problems listed in Table \cref{tab:PDE-types}. 
For brevity, define the following expression 
% common to the desirability equations for all cases of Table \cref{tab:PDE-types}
  \begin{equation*}\label{eq:L-psi}
     \mathcal{L}(\Psi)\triangleq f^{T}(\nabla_{x}\Psi)+\frac{1}{2}Tr\left(\left(\nabla_{xx}\Psi\right)\Sigma_{t}\right)
  \end{equation*}
and the function $\psi(x)$ at the boundary as  
\begin{equation*}
  \psi(x) \triangleq e^{- \frac{\phi(x)}{\lambda}} \quad x \in \partial \Omega.
\end{equation*}
% Also, define the product of state domain and time domain as $\mathcal{O} =  \Omega \times (0,T]$ and the boundary of the product domain as $\partial \mathcal{O} =  \partial \Omega \times (0,T] \cup \Omega \times \{T\}$.
% \begin{rem}

Condition \cref{eq:noise-assumption} is trivially met for systems of the form $dx_t=f(x_t)~dt +G(x_t)\left(u_t~dt\right.$ $\left.+ d\omega_t\right)$, a pervasive assumption in the adaptive control literature \cite{lavretsky2012robust} .
%Intuitively, the condition can be interpreted as the system must have sufficient control to counterbalance the effects of input noise on the system dynamics. 
This constraint restricts the design of the control penalty $R$, such that control effort
is highly penalized in subspaces with little noise, and lightly penalized in those with high
noise. Additional discussion is given in \cite[SI Sec. 2.2]{Todorov:2009wja}.

%%% ------------------------------------------------------------------------------------------- %%%
\subsection{Stochastic Control Lyapunov Functions (SCLF)}\label{sec:clf}

Before the stochastic control Lyapunov function (SCLF) is introduced, two forms of stability are defined, following the definitions in \cite[Ch. 5]{khasminskii2011stochastic}. 

\begin{definition} \label{def:stability_probability}
	Given \eqref{eq:stochastic-dynamics}, the equilibrium point at $x=0$ is stable in probability for $t\ge 0$ if for any $s\ge 0$ and $\epsilon > 0$,
	\begin{equation*}
		\lim_{x\to 0} P\left\{\sup_{t>s} |X^{x,s}(t)| > \epsilon \right\} =0
	\end{equation*}
	where $X^{x,s}$ is the trajectory of \eqref{eq:stochastic-dynamics} from $x$ at time $s$.
\end{definition}

Intuitively, \cref{def:stability_probability} is similar to the notion of stability for deterministic systems. 
%\ypl{This definition says that the sample path of \eqref{eq:stochastic-dynamics} starting from point $x$ at time $s$ will stay within a neighborhood of the origin with probability going to one as $x \to 0$. (Maybe removed if too long.)}
The following is a stronger stability definition that is similar to the notion of asymptotic stability for deterministic systems.

\begin{definition} \label{def:asymp_stability_probability}
	Given \eqref{eq:stochastic-dynamics}, the equilibrium point at $x=0$ is asymptotically stable in probability if it is stable in probability and 
	\begin{equation*}
		\lim_{x\to 0} P\left\{\lim_{t \to \infty} |X^{x,s}(t)| = 0 \right\} =1
	\end{equation*}
	where $X^{x,s}$ is the trajectory of \eqref{eq:stochastic-dynamics} from $x$ at time $s$.
\end{definition}

These notions of stability can be realized through the construction of SCLFs.
\begin{definition}  \label{def:stochastic_clf}
A stochastic control Lyapunov function for system \eqref{eq:stochastic-dynamics} is a positive
definite function $\mathcal{V} \in \mathcal{C}^{2,1}$ on a domain $\mathcal{O} =   \Omega \times \{t > 0\}$ such that 
\begin{gather*}
	\mathcal{V}(0,t) = 0, \quad \mathcal{V}(x,t)\ge \mu(|x|) \quad \forall~ t>0 \nonumber\\
	\exists~ u(x,t) \mbox{ s.t. } L(\mathcal{V}(x,t)) \leq 0  \quad \forall~(x,t) \in \mathcal{O} \backslash \{(0, t)\} \label{eq:sclf}
\end{gather*}
where $\mu \in \mathcal{K}$, and
    \begin{equation}
		L(\mathcal{V}) =\partial_t \mathcal{V}+ \nabla_x \mathcal{V}^T (f + Gu)+ \frac{1}{2}  Tr((\nabla_{xx} \mathcal{V}) B \Sigma_\epsilon B^T). \label{eq:LV}
    \end{equation}
\end{definition}
\begin{theorem} [\cite{khasminskii2011stochastic} Thm. 5.3]\label{thm:stochastic_clf}
	For system \eqref{eq:stochastic-dynamics}, assume that there exists a SCLF and a $u$ satisfying \cref{def:stochastic_clf}. Then, the equilibrium point $x = 0$ is stable in probability, and $u$ is a stabilizing controller. 
\end{theorem}

To achieve the stronger condition of asymptotic stability in probability, we have the following result.
\begin{theorem}[\cite{khasminskii2011stochastic} Thm. 5.5 and Cor. 5.1]\label{thm:stochastic_clf_asymp}
	For system \eqref{eq:stochastic-dynamics}, suppose that in addition to the existence of a SCLF and a $u$ satisfying \cref{def:stochastic_clf}, that $u$ is time-invariant, and
	\begin{gather*}
		\mathcal{V}(x,t) \le \mu'(|x|) \quad\forall~ t>0\\ 
		\quad L(\mathcal{V}(x,t)) < 0 \quad \forall~(x,t) \in \mathcal{O} \backslash \{(0, t)\}
	\end{gather*}
	where $\mu' \in \mathcal{K}$. Then, the equilibrium point $x = 0$ is asymptotically stable in probability, and $u$ is an asymptotically stabilizing controller.
\end{theorem}

% In the deterministic settings, the solution to the deterministic HJB is in fact a CLF \cite{Primbs:1999tt,freeman1996inverse}. It turns out that similar to the deterministic case, the solution to the stochastic HJB \eqref{eq:original-hjb} is a SCLF. As a result, \eqref{eq:original-hjb} emits a stabilizing controller. We will discuss this result further in Section \cref{sec:analysis}. 

%%%%%%%%%%%%%%%%%%%%%%%%%%%%%%%%%%%%%%%%%%%%%%%%%%%%%%%%%%%%%%%%%%%%%%
\subsection{Sum of Squares (SOS) Programming}\label{sec:sos}

Sum of Squares (SOS) programming is the primary tool by which approximate solutions to the HJB equation are generated in this paper. In particular, we will show how the PDE that governs the HJB may be relaxed to a set of nonnegativity constraints. SOS methods will then allow for the construction of an optimization problem where these nonnegativity constraints may be enforced. A complete introduction to SOS programming is available in \cite{Parrilo:2003fh}. Here, we review the basic definition of SOS that is used throughout the paper.
\begin{definition} \label{def:SOS}
A multivariate polynomial $f(x)$ is a SOS polynomial if there
exist polynomials $f_{0}(x),\ldots,f_{m}(x)$ such that
  \[ f(x)=\sum_{i=0}^{m}f_{i}^{2}(x). \]
  The set of SOS polynomials in $x$ is denoted as $\mathbb{S}[x]$.
\end{definition}

Accordingly, a sufficient condition for nonnegativity of a polynomial $f(x)$ is that $f(x)\in \mathbb{S}[x]$. Membership in the set $\mathbb{S}[x]$ may be tested as a convex problem \cite{Parrilo:2003fh}.  

\begin{theorem}[\cite{Parrilo:2003fh} Thm. 3.3] \label{thm:sos-test-sdp}
	The existence of a SOS decomposition of a polynomial in $n$
	variables of degree $2d$ can be decided by solving a semidefinite programming (SDP) feasibility
	problem. If the polynomial is dense (no sparsity), the dimension of the matrix inequality in the SDP is equal to $\begin{pmatrix} n+d \\ d \end{pmatrix} \times \begin{pmatrix} n+d \\ d \end{pmatrix}. $
\end{theorem}

Hence, by adding SOS constraints to the set of all positive polynomials, testing nonnegativity of a polynomial becomes a tractable SDP. The converse question, is a nonnegative polynomial necessarily a SOS, is unfortunately false, indicating that this test is conservative \cite{Parrilo:2003fh}. Nonetheless, SOS feasibility is sufficiently powerful for our purposes. \Cref{thm:sos-test-sdp} guarantees a tractable procedure to determine whether a particular polynomial, possibly parameterized, is a SOS polynomial. Our method combines multiple polynomial constraints into an optimization formulation. 
%an ability that will be granted by Schm\"{u}dgen's or Putinar's \emph{Positivstellensatz}. Before these theorems are introduced, we define 
To do so, we need to define the following polynomial sets.
\begin{definition}\label{def:preordering}
	The preordering of polynomials $g_i(x) \in \mathbb{R}[x]$ for $i = 1,2,\ldots, m$ is the set 
	\begin{equation}
		P(g_1,\ldots,g_m) = \left\{\left.\sum_{\nu \in \{0,1\}^m} s_\nu(x)g_1(x)^{\nu_1} \cdots g_m(x)^{\nu_m} \right\arrowvert s_{\nu} \in \mathbb{S}[x]\right\}.
	\end{equation}
	% The quadratic module of polynomials $g_i(x) \in \mathbb{R}[x]$ for $i = 1,2,\ldots, m$ is the set
	% \begin{equation}
	% 	M(g_1,\ldots,g_m) = \left\{\left.\sum_{i=1}^m s_i(x) g_i(x) \right\arrowvert s_i \in \mathbb{S}[x] \right\}.
	% \end{equation}
\end{definition}

The following proposition is trivial, but it is useful to incorporate the domain $\Omega$ in our optimization formulation later.
\begin{proposition}\label{lem:positive}
	Given $f(x) \in \mathbb{R}[x]$ and the domain 
	\begin{equation*}\Omega = \{x \mid  g_i(x) \in \mathbb{R}[x], g_i(x) \ge 0, i \in \{ 1,2,\ldots,m\}\},\end{equation*} 
	if $f(x) \in P(g_1,\ldots,g_m)$, then $f(x)$ is nonnegative on $\Omega$. If there exists another polynomial $f'(x)$ such that $f'(x) \ge f(x) ~\forall x \in \Omega$, then $f'(x)$ is also nonnegative on $\Omega$.
\end{proposition}
\begin{proof}
	Because $g_i(x)$ and $s_i(x)$ are nonnegative, all functions in $P(\cdot)$ are nonnegative. The second statement is trivially true given the first statement.
\end{proof}

\textbf{Example.} To illustrate an application of \cref{lem:positive}, consider a polynomial $f(x)$ defined on the domain $x\in [-1,1]$. The bounded domain can be equivalently defined by polynomials with $g_1(x) = 1+x$ and $g_2(x) = 1-x$. To certify that $f(x) \ge 0$ on the specified domain, construct a function $h(x) = s_1(x) (1+x) + s_2(x) (1-x) + s_3(x) (1+x)(1-x)$ where $s_i \in \mathbb{S}[x]$ and certify that $f(x) - h(x) \ge 0$. Notice that $h(x) \in P(1+x,1-x)$, so $h(x) \ge 0$. If $f(x) - h(x) \ge 0$, then $f(x) \ge h(x) \ge 0$. \Cref{lem:positive} is applied here. Finding the correct $s_i(x)$ is not trivial in general. Nonetheless, as mentioned earlier, if we further impose that $f(x) - h(x) \in \mathbb{S}[x]$, then checking if there exists $s_i(x)$ such that $f(x) - h(x) \in \mathbb{S}[x]$ becomes a SDP as given by \cref{thm:sos-test-sdp}. %More concretely, the procedure may begin with a limited polynomial degree for $s_i(x)$, increasing the degree until a certificate is found (if one exists) or the computation resources are exhausted. 

To simplify notation in the remainder of this text, given a domain $\Omega = \{x \mid  g_i(x) \in \mathbb{R}[x], g_i(x) \ge 0, i \in \{ 1,2,\ldots,m\}\}$, we set the notation $P(\Omega) = P(g_1,\ldots,g_m)$.% and $M(\Omega) = M(g_1,\ldots,g_m)$. 

\begin{rem}
	Depending on the computational resources available, one may choose a subset of $P(\Omega)$ to reduce the size of the resulting SDP. However, the chances of finding a certificate reduces as a consequent. This polynomial set is often used in the discussions of Schm\"{u}dgen's Positivstellensatz, which states that if $f(x)$ is positive on a compact domain $\Omega$, then $f(x) \in P(\Omega)$ \cite{lasserre2008nonlinear,Parrilo:2003fh}.
	% Choosing either $M(\Omega)$ or $P(\Omega)$ relies on the computational resources available. Although $M(\Omega) \subset P(\Omega)$ and therefore the chances of finding a certificate is larger using $P(\Omega)$, the resulting SDP is also larger. In addition, using other subsets of $P(\Omega)$ apart from $M(\Omega)$ does not change the results. These polynomial sets are often used in the discussions of Schm\"{u}dgen's or Putinar's Positivstellensatz. Loosely speaking, Schm\"{u}dgen's Positivstellensatz states that if $f(x)$ is positive on a compact domain $\Omega$, then $f(x) \in P(\Omega)$ \cite{lasserre2008nonlinear,Parrilo:2003fh}. 
	%These theorems may fail in our case because the function of interest is nonnegative instead of positive. A more general Stengle's {Positivstellensatz} may fail as well because of the same reason. 
% \ypl{But, we don't really care because we are not trying to prove that $f(x)$ is positive or nonnegative. We are trying to find $f(x)$ and $s_i(x)$ such that $f(x) - h(x) \in \mathbb{S}[x]$. Since $f(x)$ is a decision variable, there is a decent chance that $f(x) - h(x) \in \mathbb{S}[x]$ I think. Not sure how to say all these rigorously yet.}	
	%our function of interest is also an optimization variable in the semidefinite program. Section \ypl{XXXXX} will further explain this point after the optimization problem is introduced.
\end{rem}

%%% ------------------------------------------------------------------------------------------- %%%
\section{SOS Relaxation of the HJB PDE}\label{sec:problem}

SOS programming has found many uses in combinatorial optimization, control theory, and
other applications. This section now adds solving the linear HJB to this list.
We would like to emphasize the following standing assumption, necessary in moment and SOS-based methods \cite{lasserre2008nonlinear,Parrilo:2003fh}.
\begin{ass} \label{ass:domain}
	Assume that system \eqref{eq:stochastic-dynamics} evolves on a compact domain $\Omega \subset \mathbb{R}^n$, and $\Omega$ is a basic closed semialgebraic set such that $\Omega = \{x \mid g_i(x) \in \mathbb{R}[x], g_i(x) \ge 0, i \in \{1, \ldots, k\}\}$ for some $k \ge 1$. Then, the boundary $\partial\Omega$ is polynomial representable. We use the notation $\partial\Omega = \{x \mid h_i(x) \in \mathbb{R}[x], \prod_{i=1}^m h_i(x)=0\}$ for some $m \ge 1$ to describe the boundary. 
\end{ass}

% For ease of notation, we also define the operator $\mathcal{D}$ ($\mathcal{B}$) that maps a semialgebraic domain (boundary) into a polynomial. This will be useful in the application of Lemma \cref{lem:positive}. 
The following definitions formalize several operators that are useful in the sequel.
\begin{definition}\label{def:domain}
	Given a basic closed semialgebraic set $\Omega = \{x \mid g_i(x) \in \mathbb{R}[x], $ $g_i(x) \ge 0, i \in \{1, \ldots, k\}\}$ and a set of SOS polynomials,
	\begin{equation*}\mathcal{S} = \{s_\nu(x)\mid s_\nu(x) \in\mathbb{S}[x], \nu \in \{0,1\}^k\},\end{equation*} 
	define the operator $\mathcal{D}$ as 
	\begin{equation*}
		\mathcal{D}(\Omega,\mathcal{S}) = \sum_{\nu \in \{0,1\}^k} s_\nu(x) g_1(x)^{\nu_1}\cdots g_k(x)^{\nu_k}
	\end{equation*}
	where $s_\nu \in \mathcal{S}$ and $\mathcal{D}(\Omega,\mathcal{S}) \in P(\Omega)$.
\end{definition}

\begin{definition}\label{def:boundary}
	Given a polynomial inequality, $p(x) \ge 0$ defined on $\Omega$, the boundary of a compact set $\partial\Omega = \{x \mid h_i(x) \in \mathbb{R}[x], \prod_{i=1}^m h_i(x)=0\}$ and a set of  polynomials, 
	\begin{equation*}\mathcal{T} = \{t_i(x) \mid t_i(x) \in\mathbb{R}[x], i \in \{1, \ldots, m\}\},\end{equation*}  
	define the operator $\mathcal{B}$ as 
	\begin{equation*}
		\mathcal{B}(p(x),\partial \Omega,\mathcal{T}) = \{p(x)- t_{i}(x) h_{i}(x) \mid  i \in \{1, \ldots, m\}\}
	\end{equation*}
	where $t_i \in \mathcal{T}$ and $\mathcal{B}$ returns a set of polynomials that is nonnegative on $\partial \Omega$.
\end{definition}

%%%%%%%%%%%%%%%%%%%%%%%%%%%%%%%%%%%%%%%%%%%%%%%%%%%%%%%%%%%%%%%%%%%%%%
\subsection{Relaxation of the HJB equation}\label{sec:hjb-relax}
If the linear HJB \eqref{eq:hjb-pde} is not uniformly parabolic \cite{crandall1992user}, a classical solution may not exist. The notion of \emph{viscosity solutions} is developed to generalize the classical solution. We refer readers to \cite{crandall1992user} for a general discussion on viscosity solutions and \cite{Fleming:2006tl} for a discussion on viscosity solutions related to Markov diffusion processes. %Due to the depth of the topic, and our direct use of existing results, we direct the reader to these references for a review. Nonetheless, the definition of viscosity solution is reviewed next.
 % \ypl{Add in the definition again because Joel said so.}

\begin{definition} [\cite{crandall1992user} Def. 2.2]\label{def:viscosity solutions}
	Given $\Omega \subset \mathbb{R}^N$ and a partial differential equation
\begin{equation}\label{eq:general_pde}
F(x,u,\nabla_x u,\nabla_{xx} u) = 0 
\end{equation}
where $F:\mathbb{R}^N \times \mathbb{R} \times \mathbb{R}^N \times \mathcal{S} (N) \to \mathbb{R}$, $\mathcal{S}(N)$ is the set of real symmetric $N\times N$ matrices, and $F$ satisfies
\begin{equation*}
 	   F(x,r,p,X) \le F(x,s,p,Y) \mbox{ whenever $r\le s$ and $Y\le X$},
\end{equation*}
then a \textbf{viscosity subsolution} of \eqref{eq:general_pde} on $\Omega$ is a function  
% $u \in \mathcal{C}^2(\Omega)$
$u \in USC(\Omega)$ such that
\begin{equation*}\label{eq:subsolution}
F(x,u,\nabla_x u,\nabla_{xx} u) \le 0\quad \forall ~x\in\Omega, (p,X) \in J^{2,+}_\Omega u(x).
\end{equation*}
Similarly, a \textbf{viscosity supersolution} of \eqref{eq:general_pde} on $\Omega$ is a function 
% $u \in \mathcal{C}^2(\Omega)$
 $u\in LSC(\Omega)$ 
such that
\begin{equation*}\label{eq:supersolution}
F(x,u,\nabla_x u,\nabla_{xx} u) \ge 0 \quad \forall~ x\in\Omega, (p,X) \in J^{2,-}_\Omega u(x).
\end{equation*}
Finally, $u$ is a \textbf{viscosity solution} of \eqref{eq:general_pde} on $\Omega$ if it is both a viscosity subsolution and a viscosity supersolution in $\Omega$.
\end{definition}
%The notation $\mathcal{C}^2(\Omega)$ indicates the function is twice differentiable on the domain $\Omega$.

The notations $USC(\Omega)$ and $LSC(\Omega)$ represent the sets of upper and lower semicontinuous functions on domain $\Omega$ respectively, and $J^{2,+}_\Omega u(x)$ and $J^{2,-}_\Omega u(x)$ represents the second order ``superjets'' and ``subjets'' of $u$ at $x$ respectively, a completely unrestrictive domain in our setting. For further details, readers may refer to \cite{crandall1992user}. 
% Note that if $u$ is a smooth function such that $F(x,u,\nabla_x u,\nabla_{xx} u) \le 0$, then $u$ satisfies the definition of viscosity subsolution. Similarly, if $u$ is a smooth function such that $F(x,u,\nabla_x u,\nabla_{xx} u) \ge 0$, then $u$ satisfies the definition of viscosity supersolution. 
For the remainder of this paper, we assume a unique nontrivial viscosity solution to \eqref{eq:hjb-pde-value} and \eqref{eq:hjb-pde} exists (see \cite{Fleming:2006tl}, Chapter V) and denote them as $V^*$ and $\Psi^*$ respectively.

% \mbh{I think we can get away without talking about uniqueness of the true solution. If you think it should be included, I have commented out a theorem that may help.}
% \begin{theorem}
% There exists a unique viscosity solution $\Psi^*$ to \eqref{eq:hjb-pde}.
% \end{theorem}
% \begin{proof}
% This result is a trivial application of \cite{crandall1992user}, Theorem 7.9
% \end{proof}

%  We do not include additional technical assumptions that ensure the existence of a unique viscosity solution because it is not the main focus of this paper. Readers may refer to \cite{crandall1992user,Fleming:2006tl} for details. Henceforth, we only consider \eqref{eq:hjb-pde} that has a unique nontrivial viscosity solution $\Psi^*$, and we refer to the viscosity solution $\Psi^*$ as the solution to \eqref{eq:hjb-pde}. This assumption also implies that \eqref{eq:hjb-pde-value} has a nontrivial viscosity solution because the logaritheoremic transformation is bijective.
The equality constraints of \eqref{eq:hjb-pde} may be relaxed as follows
\begin{subequations}
  \begin{gather}
    \frac{1}{\lambda}q\Psi - \mathcal{L}(\Psi) \le (\ge) 0 \label{eq:over-approximation-relax-1} \\
		\Psi(x) \le (\ge) \psi(x) \qquad  x \in \partial \Omega. 
  \end{gather}\label{eq:over-approximation-relax}
\end{subequations}
Such a relaxation provides a point-wise bound to the solution $\Psi^*$, and this relaxation may be enforced via SOS programming. In particular, a solution to \eqref{eq:over-approximation-relax}, denoted as $\Psi_l (\Psi_u)$, is a lower (upper) bound on the solution $\Psi^*$ over the entire problem domain. 

\begin{theorem} \label{prop:Psi_bound}
Given a smooth function $ \Psi_l (\Psi_u)$ that satisfies \eqref{eq:over-approximation-relax}, then $\Psi_l (\Psi_u)$ is a viscosity subsolution (supersolution) and $\Psi_l \le \Psi^{*} (\Psi_u \ge \Psi^{*})$ for all $x \in \Omega$.
\end{theorem}

\begin{proof}
	By \cref{def:viscosity solutions}, the solution $\Psi_l$ is a viscosity subsolution where $F$ in \eqref{eq:general_pde} is given by \eqref{eq:over-approximation-relax-1}. Note that $\Psi^*$ is both a viscosity subsolution and a viscosity supersolution, and $\Psi_l \le \Psi^{*}$ on the boundary $\partial \Omega$. Hence, by the maximum principle \cite[Thm. 3.3]{crandall1992user}, $\Psi_l \le \Psi^{*}$ for all $x \in \Omega$. The proof is identical for $\Psi_u$.
%According to Theorem \cref{thm:viscosity_bound}, the solution $\Psi_l$ is a viscosity subsolution by construction whereby $F(t,u,\nabla_x u,\nabla_{xx} u) =  - \partial_t u + \frac{1}{\lambda} q u - \mathcal{L}(u) $. The second order term $\Sigma_\epsilon$ is always positive-semidefinite as it is a covariance, ensuring the PDE is degenerate elliptic. Furthermore, the positivity of the state cost ensures the system is proper. The only other requirement for the desired guarantee is that for $\Psi_l$ to be upper-semicontinuous. This requirement is satisfied as the solution is in fact a polynomial and therefore is infinitely differentiable.
\end{proof}
%
% Similarly, the analogous relaxation
%   \begin{gather}
%      \frac{1}{\lambda}q\Psi -\mathcal{L}(\Psi) \ge 0\label{eq:over-approximation-relax-upper}\\
% 		\Psi(x) \ge  \psi(x) \qquad (x) \in \partial \Omega \nonumber
%   \end{gather}
% gives an over-approximation of the desirability function, and its solution, denoted as $\Psi_u$,
% is an upper bound of $\Psi^*$ over domain $\Omega$. Thus, we also have
%
% \begin{theorem}  \label{prop:Psi_bound_upper}
% Given a smooth function $ \Psi_u $ that satisfies \eqref{eq:over-approximation-relax-upper}, then $\Psi_u$ is a viscosity supersolution and $\Psi_u \ge \Psi^{*}$ for all $x \in \Omega$.
% \end{theorem}
% \begin{proof}
% 	The proof is identical to the proof of \cref{prop:Psi_bound}.
% \end{proof}

Because the logarithmic transform \eqref{eq:log-transform} is monotonic, one can relate
these bounds on the desirability function to bounds on the value function as follows
\begin{proposition} \label{prop:v-upper-lower}
If the solution to \eqref{eq:hjb-pde-value} is $V^*$, given solutions $V_u = -{\lambda}\log \Psi_l$ and $V_l = -{\lambda}\log \Psi_u$ from \eqref{eq:over-approximation-relax}, then $V_u \ge V^*$ and $V_l \le V^*$.
\end{proposition}
\begin{proof}
	Recall that $V^* = -{\lambda}\log \Psi^*$. Apply \cref{prop:Psi_bound}, $V_u \ge V^*$ and $V_l \le V^*$.
\end{proof}

The solutions to \eqref{eq:over-approximation-relax} do not satisfy \eqref{eq:hjb-pde} exactly, but they provide point-wise bounds to the solution $\Psi^*$. %Indeed, these super or sub-solutions are \emph{viscosity} solutions \cite{crandall1992user, Fleming:2006tl}, and have been studied primarily as a vehicle to generalize the notion of a solution to a partial differential equation. %The procedure outlined above provides a method to generate these viscosity solutions, which we later show are of practical interest.

% \begin{rem}
	% A viscosity solution is always continuous but not necessarily smooth, and not all viscosity subsolutions and supersolutions are continuous and smooth. However, the viscosity subsolutions and supersolutions that are constructed in our optimization scheme in the next section are polynomials. And, thus, they are always smooth and continuous within the problem domain with well defined classical derivatives. 
	%For a more general approach that considers discontinuous viscosity solutions and discontinuous Lyapunov functions, refer to \cite{bardi2005almost}. \ypl{Do we need to mention this last part?}
	%\mbh{I don't think so}
% \end{rem}

%%%%%%%%%%%%%%%%%%%%%%%%%%%%%%%%%%%%%%%%%%%%%%%%%%%%%%%%%%%%%%%%%%%%%%
\subsection{SOS Program}\label{sec:oc_synthesis}

Given that relaxation \eqref{eq:over-approximation-relax} results in a point-wise upper and lower bound to the exact solution of \eqref{eq:hjb-pde}, we construct the following optimization problem that provides a suboptimal controller with bounded residual error:
%
% While we show below that the approximating desirability functions converge to the optimal
% desirability functions with increasing polynomial degree, one must consider the closed loop system
% behavior when using a suboptimal viscosity solution.  The optimal control is based on the gradient
% of the value function (see Equation (\cref{eq:optimal-u})), rather than the value function itself.
% In general, there is no guarantee on the convergence of the gradient of an approximating function to
% the gradient of the optimal value function.  Nonetheless, we show below that a weak convergence of
% the optimal cost can be guaranteed if the relaxation of \eqref{eq:hjb-pde} is recast into the
% following form
  \begin{alignat}{2}
	\min_{\Psi_l,\Psi_u} \quad& \epsilon \label{eq:hjbjoin} \\
	s.t. \quad &  \frac{1}{\lambda} q \Psi_l 
              - \mathcal{L}(\Psi_l) \leq 0  && \quad x \in\Omega\nonumber \\
 	& 0 \leq  \frac{1}{\lambda} q \Psi_u 
              - \mathcal{L}(\Psi_u) &&  \quad x \in \Omega \nonumber \\
	& \Psi_u - \Psi_l \leq \epsilon &&\quad x\in \Omega \nonumber \\
	%& Tr\left(\left(\nabla_{xx}\Psi_l\right)\Sigma_{t}\right) \leq 0  && \quad x \in\Omega \nonumber \\
	& 0 \le \Psi_l  \le\psi  \le \Psi_u   && \quad x \in\partial \Omega \nonumber \\
	% & \Psi_l \le \rho(x) &&  \quad x \in \Omega \nonumber
	& \partial_{x^i} \Psi_l \le 0 && \quad x^i \ge 0 \nonumber \\
	& \partial_{x^i} \Psi_l \ge 0 && \quad x^i \le 0 \nonumber \\
	& \Psi_l(0) = 1  \nonumber 
  \end{alignat} 
  where $x^i$ is the $i$-th component of $x \in \Omega$. As mentioned in \cref{sec:hjb-relax}, the first two constraints result from the relaxations of the HJB equation, and the fourth constraint
arises from the relaxation of the boundary conditions. The third constraint ensures that the difference between the upper bound and lower bound
solution is bounded, and the last three constraints ensure that the solution yields a stabilizing controller, as will be made clear in \cref{sec:analysis}. Note that in the optimization problem, $\Psi_u$ and $\Psi_l$ are polynomials whereby the coefficients and the degree for both are optimization variables. The term $\epsilon$ is related to the error of the approximation. 

As discussed in the review of SOS techniques, a general optimization problem involving parameterized nonnegative polynomials is not necessarily tractable. In order to solve \eqref{eq:hjbjoin} using a polynomial-time algorithm, we restrict the polynomial inequalities such that they are SOS polynomials instead of nonnegative polynomials. We therefore apply \cref{lem:positive} to relax optimization problem \eqref{eq:hjbjoin} into
  \begin{alignat}{2}
	\min_{\Psi_l,\Psi_u,\mathcal{S},\mathcal{T}} \quad& \epsilon \label{eq:hjbjoin-sos} \\
	s.t. \quad & - \frac{1}{\lambda} q \Psi_l 
              + \mathcal{L}(\Psi_l)- \mathcal{D}(\Omega,\mathcal{S}_1) \in \mathbb{S}[x]\nonumber \\
 	&  \frac{1}{\lambda} q \Psi_u 
              - \mathcal{L}(\Psi_u) - \mathcal{D}(\Omega,\mathcal{S}_2) \in \mathbb{S}[x] \nonumber \\
	& \epsilon - (\Psi_u - \Psi_l ) - \mathcal{D}(\Omega,\mathcal{S}_3)  \in \mathbb{S}[x]\nonumber \\
	%& -Tr\left(\left(\nabla_{xx}\Psi_l\right)\Sigma_{t}\right) - \mathcal{D}(\Omega,\mathcal{S}_4) \in \mathbb{S}[x] \nonumber \\	
%	& 0\le \Psi_l  \le e^{-\frac{\phi_T(x_T)}{\lambda}}  \le \Psi_u \quad x\in\partial\Omega \nonumber
	&  \mathcal{B}( \Psi_l ,\partial\Omega,\mathcal{T}_1) \in \mathbb{S}[x]\nonumber \\
	&  \mathcal{B}(\psi -\Psi_l ,\partial\Omega,\mathcal{T}_2)\in \mathbb{S}[x] \nonumber \\
	&  \mathcal{B}( \Psi_u - \psi,\partial\Omega,\mathcal{T}_3) \in \mathbb{S}[x]\nonumber \\
	& -\partial_{x^i} \Psi_l - \mathcal{D}(\Omega \cap \{x^i \ge 0\}, \mathcal{S}_4) \in \mathbb{S}[x]\nonumber \\
	& \partial_{x^i} \Psi_l  - \mathcal{D}(\Omega \cap \{-x^i \ge 0\}, \mathcal{S}_5) \in \mathbb{S}[x]\nonumber \\
	& \Psi_l(0) = 1  \nonumber 
  \end{alignat}
  where $\mathcal{S} = (\mathcal{S}_1,\ldots,\mathcal{S}_4,\mathcal{S}_5)$, $\mathcal{S}_i \subseteq \mathbb{S}[x]$ is defined as in \cref{def:domain}, $\mathcal{T} = (\mathcal{T}_1,\mathcal{T}_2,\mathcal{T}_3)$, and $\mathcal{T}_j \subseteq \mathbb{R}[x]$ is defined as in \cref{def:boundary}. With a slight abuse of notation, $\mathcal{B}(\cdot) \in \mathbb{S}[x]$ implies that each polynomial in $\mathcal{B}(\cdot)$ is a SOS polynomial.

  % In the above optimization, the optimization variables are the coefficients of $\Psi_u$, $\Psi_l$, and polynomials in $\mathcal{S}$ and $\mathcal{T}$, while the degree of these polynomials is a meta-optimization variable chosen by the practitioner.  In the most general form, $\mathcal{S}_i $ are four different sets of SOS polynomials, and $\mathcal{T}_j $ are three different sets of polynomials. However, to limit the size of optimization variables, one can choose to set the SOS polynomial set $\mathcal{S}_i $ to be the same for all $i$. Similarly, one can also set the polynomial set $\mathcal{T}_j$ to be the same for all $j$. 
  
% At this point, the optimization problem \eqref{eq:hjbjoin-sos} is ``almost'' a tractable problem.
 If the polynomial degrees are fixed, optimization problem \eqref{eq:hjbjoin-sos} is convex and solvable using a semidefinite program via \cref{thm:sos-test-sdp}. The next section will discuss the systematic approach we used to solve the optimization problem. Henceforth, denote the solution to \eqref{eq:hjbjoin-sos} as $(\Psi_u,\Psi_l, \mathcal{S},\mathcal{T}, \epsilon)$. 
 % In addition, when the optimization problem \eqref{eq:hjbjoin-sos} is referred to in later text, we mean the optimization problem when the degrees of polynomials are fixed, and the optimization variables are the coefficients of $\Psi_u$, $\Psi_l$, and polynomials in $\mathcal{S}$ and $\mathcal{T}$. 

% \subsection*{Existence of optimization solutions}
\begin{rem} \label{rem:continuous}
By \cref{def:viscosity solutions}, the viscosity solution is a continuous function. Consequently, the solution $\Psi^*$ is a continuous function defined on a bounded domain. Hence, $\Psi_u$ and $\Psi_l$ can be made arbitrary close to $\Psi^*$ by the Stone-Weierstrass Theorem \cite{rudin1964principles} in \eqref{eq:hjbjoin}. However, this guarantee is lost when $\Psi_u$ and $\Psi_l$ are restricted to be a SOS polynomials. The feasible set of the optimization problem \eqref{eq:hjbjoin-sos} is therefore not necessarily non-empty for a given polynomial degree. One would not expect feasibility for all instances of \eqref{eq:hjbjoin-sos} as this would imply there exists is a linear stabilizing controller for any given system. %\mbh{Tried to give a simpler explanation.}
% \ypl{Joel: Can we say the rest more simply?} This infeasibility for certain degree is expected. If the feasible set is non-empty for all given degrees, then there exists an upper bound $V^u$ that yields a stabilizing controller for any polynomial degree, an exceedingly strong statement.
\end{rem}

% The feasible set of the optimization problem \eqref{eq:hjbjoin-sos} is not necessarily non-empty for a given polynomial degree. This infeasibility for certain degree is expected because if the feasible set is non-empty for all given degree, then there exists an upper bound $V^u$ that emits a stabilizing controller for any polynomial degree. Such strong statement cannot be true for arbitrary nonlinear systems. Determining the minimum degree required for the optimization problem \eqref{eq:hjbjoin-sos} to be feasible will be explored in the future. 

%%%%%%%%%%%%%%%%%%%%%%%%%%%%%%%%%%%%%%%%%%%%%%%%%%%%%%%%%%%%%%%%%%%%%%
\subsection{Controller Synthesis}\label{sec:controller_synthesis}

Let $d$ be the maximum degree of $\Psi_l$, $\Psi_u$ and polynomials in $\mathcal{S}$ and $\mathcal{T}$, and denote $(\Psi^d_u,\Psi^d_l, \mathcal{S}^d,\mathcal{T}^d, \epsilon^d)$ as a solution to \eqref{eq:hjbjoin-sos} when the maximum polynomial degree is fixed at $d$. The hierarchy of SOS programs with increasing polynomial degree produces a sequence of (possibly empty) solutions $(\Psi^d_u,\Psi^d_l, \mathcal{S}^d,\mathcal{T}^d, \epsilon^d)_{d\in I}$, where $I \subset \mathbb{Z}_+$. This sequence will be shown in the next section to improve, under the metric of the objective in \eqref{eq:hjbjoin-sos}.

In other words, if solutions exist for $d$ and $d'$ such that $d > d'$, then $\epsilon^d \le \epsilon^{d'}$. Therefore, one could keep increasing the degree of polynomials in order to achieve tighter bounds on $\Psi^*$, and invariably, $V^*$. The use of such hierarchies has become commonplace in polynomial optimization \cite{lasserre2001global,Parrilo:2003fh}. If at certain degree, $\epsilon^d = 0$, the solution $\Psi^*$ is found.  %The next section will analyze the performance of this optimization hierarchy.

Once a satisfactory error is achieved or computational resources run out, the lower bound $\Psi^d_l$ can be used to compute a suboptimal controller where $d$ is the maximum degree computed. Recall that $u^* = -R^{-1}G^T \nabla_x V^*$ and $V^* = -{\lambda} \log \Psi^*$. The suboptimal controller $u^\epsilon$ for a given degree $d$ and error $\epsilon^d$ is computed as $u^{\epsilon^d} = -R^{-1}G^T \nabla_x V^d_u$ where $V^d_u = - {\lambda} \log {\Psi^d_l}$. Even when $\epsilon^d$ is larger than a desired value, the solution $\Psi^d_l$ still satisfies conditions in \cref{def:stochastic_clf} to yield a stabilizing suboptimal controller. Next section will analyze properties of the solutions and the suboptimal controller.
\section{Analysis}\label{sec:analysis}

This section establishes several properties of the solutions to the optimization problem \eqref{eq:hjbjoin-sos} that are useful for feedback control. First we show that the solutions in the SOS program hierarchy are uniformly bounded relative to the exact solutions. We next prove that the relaxed solutions to the stochastic HJB equation are SCLFs, and the approximated solution leads to a stabilizing controller. Finally, we show that the costs of using the approximate solutions as controllers are bounded above by the approximated value functions. 

%%%%%%%%%%%%%%%%%%%%%%%%%%%%%%%%%%%%%%%%%%%%%%%%%%%%%%%%%%%%%%%%%%%%%%
\subsection{Properties of Approximated Desirability Functions} \label{sec:prop-desirability}

First, the approximation error of $\Psi_l$ or $\Psi_u$ obtained from \eqref{eq:hjbjoin-sos} is computed relative to the true desirability function $\Psi^*$.
\begin{proposition} \label{prop:psi_estimate}
Given a solution $(\Psi^d_u,\Psi^d_l, \mathcal{S}^d,\mathcal{T}^d,\epsilon^d)$ to \eqref{eq:hjbjoin-sos} for a given degree $d$, the approximation
error of the desirability function is bounded as $||\Psi^d - \Psi^*||_{\infty} \leq \epsilon^d$ where
$\Psi^d$ is either $\Psi^d_u$ or $\Psi^d_l$.
\end{proposition}

\begin{proof}
By \cref{prop:Psi_bound}, $\Psi^d_l$ is the lower bound of $\Psi^*$,
and $\Psi^d_u$ is the upper bound of $\Psi^*$. So, $\epsilon^d \geq \Psi^d_u - \Psi^d_l \geq 0$ and $\Psi^d_u
\geq \Psi^* \geq \Psi^d_l$. Combining both inequalities, one has $\Psi^d_u - \Psi^* \leq \epsilon^d$
and $\Psi^*- \Psi^d_l \leq \epsilon^d$. Therefore, $||\Psi^d - \Psi^*||_{\infty} \leq \epsilon^d$ where
$\Psi^d$ is either $\Psi^d_u$ or $\Psi^d_l$.
\end{proof}

\begin{proposition} \label{prop:gap}
%Let $\Psi^d_l$ and $\Psi^d_u$ be polynomial approximations of the desirability function with maximum
%polynomial degree $d$. 
The hierarchy of SOS programs consisting of solutions to \eqref{eq:hjbjoin-sos} with
increasing polynomial degree produces a sequence of solutions $(\Psi^d_u,\Psi^d_l, \mathcal{S}^d,\mathcal{T}^d, \epsilon^d)$ such that $\epsilon^{d+1} \le \epsilon^{d}$ for all $d$.
\end{proposition}

% \jwb{This next proof needs a lot of work? How do I know that ``$\gamma_u^{i} +
% \gamma_l^{i} + \epsilon^i$, decreases monotonically''?  This must be proven, as it is not obvious. Also,
% the symbols $\epsilon_u^i$ and $\epsilon_u^j$ are not defined, even though their definitions are obvious.}

\begin{proof}
	Polynomials of degree $d$ form a subset of polynomials of degree $d+1$. Thus, at a higher polynomial degree $d+1$, a previous solution at a lower polynomial degree $d$ is still a feasible solution when the coefficients for monomials with total degree $d+1$ is set to 0. Consequently, the optimal value $\epsilon^{d+1}$ cannot be larger than $\epsilon^{d}$ for all $d$. 
% 	Thus, if $i > $
% With increasing polynomial degree of the approximation, the objective function, $\gamma_u^{i} +
% \gamma_l^{i} + \epsilon^i$, decreases monotonically. Let $s^i = \gamma_u^{i} + \gamma_l^{i} +
% \epsilon^i$ and $s^j = \gamma_u^{j} + \gamma_l^{j} + \epsilon^j$ where $i > j$. Suppose $i$ is
% increased from $j$ until $s^i \leq \epsilon_u^j$. Then, $\epsilon_u^i \leq \epsilon_u^j$ because
% $\epsilon_u^i \leq s^i$. Similar logic applies to both $\epsilon_l$ and $\gamma$.  Thus,
% $\gamma_u^{i}$, $\gamma_l^{i}$ and $\epsilon^i$ decreases eventually as $i$
% increases.
\end{proof}

Thus, as the polynomial degree of the optimization problem is increased, the pointwise error $\epsilon$ is non-increasing. Therefore, one could keep increasing the degree of polynomials in order to achieve tighter bounds on $\Psi^*$, and invariably, $V^*$. However, $\epsilon$ is only non-increasing as the polynomial degree is increased, and a convergence of the bound $\epsilon$ to zero is not guaranteed because we restrict the approximating space to SOS. The possible lack of convergence to zero is the trade off for an efficient algorithm.

Although the \emph{bound} on the pointwise error is non-increasing, the actual difference between $\Psi$ and $\Psi^*$ may increase between iterations. %We bound this variation as follows.
\begin{corollary} \label{col:gap}
Suppose $||\Psi^{d} - \Psi^*||_{\infty} \le \epsilon^{d}$ and $||\Psi^{d+1} - \Psi^*||_{\infty} = \gamma^{d+1}$. Then, $\gamma^{d+1} \le \epsilon^{d}$. 
\end{corollary}

\begin{proof}
By \cref{prop:gap}, $\epsilon ^{d+1} \le \epsilon^d$. Because $\gamma^{d+1} \le \epsilon^{d+1}$, $\gamma^{d+1} \le \epsilon^{d}$
\end{proof}

In other words, the approximation error of the desirability function for a SOS program using $d+1$ polynomial degree cannot increase such that it is larger than $\epsilon^{d}$ in each step of the hierarchy of SOS programs, which is non-increasing.

% In addition, notice that the third constraint in \eqref{eq:hjbjoin-sos} plays an important role to ensure the non-increasing properties of the solutions. Therefore, it is necessary to include this constraint in the optimization.

%%%%%%%%%%%%%%%%%%%%%%%%%%%%%%%%%%%%%%%%%%%%%%%%%%%%%%%%%%%%%%%%%%%%%%
\subsection{Properties of Approximated Value Functions}\label{sec:cost_convergence}

Up to this point, the analysis has focused on properties of the desirability solution. We now investigate the implications of these results upon the value function, which is related to the desirability via the logarithmic transform \eqref{eq:log-transform}. Henceforth, denote the solution to \eqref{eq:hjb-pde-value} as $V^*(x_t) = \min_{u[t:T]}\mathbb{E}_{\omega_t}[J(x_t)] = - \lambda \log \Psi^*(x_t)$, the solution to \eqref{eq:hjbjoin-sos} for a fixed degree $d$ as $(\Psi_u,\Psi_l, \mathcal{S},\mathcal{T}, \epsilon)$, and the suboptimal value function computed from the solution of \eqref{eq:hjbjoin-sos} as $V_u = - {\lambda} \log {\Psi_l}$. Only $\Psi_l$ and $V_u$ are considered henceforth, because $\Psi_l$, but not $\Psi_u$, gives an approximate value function that satisfies the properties of SCLF in \cref{def:stochastic_clf}, a fact shown in the next section.

\begin{theorem} \label{thm:cost-upper}
For all $x\in \Omega$, $V_u$ is an upper bound of $V^*$ such that
  \begin{equation}
	 0 \leq V_u - V^* \leq -\lambda \log\left(1-\min\left\{1,\frac{\epsilon}{\eta}\right\}\right) 
  \end{equation}
where $\eta = e^{-\frac{\norm{V^*}_\infty}{\lambda}}$.
\end{theorem}
\begin{proof}
	% First, show that $V_u - V^* \ge 0$.
% 	By fourth constraint in \eqref{eq:hjbjoin-sos}, $\Psi_l \ge 0$. If $\Psi_l(0)>0$, the fraction $\frac{\Psi_l}{\Psi_l(0)} \ge 0$, and it satisfies \eqref{eq:over-approximation-relax} because \eqref{eq:over-approximation-relax} is linear with respect to $\Psi$. \ypl{What to do about the boundary condition? if the boundary condition is not satisfy, the next sentence is not true. Maybe can't do scaling after all? Or maybe can figure out a way to scale the boundary?} Then, by Theorem \cref{prop:Psi_bound}, $\frac{\Psi_l}{\Psi_l(0)}  \le \Psi^{*}$. Through the logaritheoremic transformation, 
By \cref{prop:v-upper-lower}, $V_u \ge V^* $ and hence, $V_u - V^* \ge 0$. To prove the other inequality, by \cref{prop:psi_estimate},
  \begin{align*}
	V_u - V^*  &= - \lambda \log \frac{\Psi_l}{\Psi^*} \leq -\lambda \log \frac{\Psi^*-\epsilon}{\Psi^*} \leq -\lambda \log \left(1-\frac{\epsilon}{\eta}\right) .
  \end{align*}	
The last inequality holds because $\Psi^* \geq e^{-\frac{\norm{V^*}_\infty}{\lambda}}$ by definition
in \eqref{eq:log-transform}. Since $\Psi_l$ is the lower bound of $\Psi^*$, the right hand side of the
first equality is always a positive number. Therefore, $V_u$ is a point-wise upper bound of
$V^*$. 
% \mbh{I think we can skip the next sentence, it invites more questions than it answers. For $\epsilon>\eta$ it's just a useless bound.}
% When $\epsilon > \eta$, we define the bound to be positive infinity. \ypl{Joel: Need to say when does this happen. Discuss more.}
\end{proof}

% \ypl{If the boundary of $V$ is set at infinity, $\log\left(1-\min\left\{1,\frac{\epsilon}{\eta}\right\}\right)$ will always be $- \infty$. So, the upper bound is meaningless.}

\begin{corollary} \label{col:cost-upper}
%Given a solution $(\Psi_u,\Psi_l, \mathcal{S},\mathcal{T},\epsilon)$ to \eqref{eq:hjbjoin-sos}, 
Let $V_u^{d} = -\lambda \log \Psi^{d}_l$ and $V_u^{d+1} = -\lambda \log \Psi^{d+1}_l$. If $\Psi_u^d - \Psi^*\le \epsilon^d$ and $V_u^{d+1} - V^* = \gamma^{d+1}$, then $\gamma^{d+1} \le -\lambda \log\left(1-\min\left\{1,\frac{\epsilon^{d}}{\eta}\right\}\right)$.
%is non-increasing as the degree of polynomial $d$ increases in the hierarchy of SOS programs.
\end{corollary}

\begin{proof}
% If the data in the HJB equation, the PDE coefficients and boundary conditions, are smooth, as has
% been assumed, then the associated viscosity solutions are also guaranteed to be smooth, and
% therefore may be approximated by polynomials \cite{crandall1992user}. 
This result is given by \cref{col:gap,thm:cost-upper}.
\end{proof}

At this point, we have shown that the lower bound of the desirability function yields an upper
bound of the suboptimal cost. More importantly, the upper bound of the suboptimal cost is not increasing as the degree of polynomial increases. 

%%%%%%%%%%%%%%%%%%%%%%%%%%%%%%%%%%%%%%%%%%%%%%%%%%%%%%%%%%%%%%%%%%%%%%
\subsection{Approximate HJB solutions are SCLFs} \label{sec:hjb is sclf}

This section shows that the approximate value function derived from the approximation,
$\Psi_l$, is a SCLF.

\begin{theorem}\label{thm:lower-sclf}
$V_u$ is a stochastic control Lyapunov function according to \cref{def:stochastic_clf}. 
\end{theorem}

\begin{proof}
	The constraint $\Psi_l(0) = 1$ in \eqref{eq:hjbjoin-sos} ensures that $V_u(0) = - \lambda \log {\Psi_l(0)} = 0$. Notice that all terms in $J(x,u)$ from \eqref{eq:cost-functional} are positive definite, resulting in $V^*$ being a positive definite function. In addition, by \cref{prop:v-upper-lower}, $V^u \ge V^*$. Hence, $V^u$ is also a positive definite function. The second and third to last constraints in \eqref{eq:hjbjoin-sos} ensures that $\Psi_l$ is nonincreasing away from the origin. Hence, $V_u$ is nondecreasing away form the origin satisfying $\mu(|x|) \le V_u(x) \le \mu'(|x|)$ for some $\mu, \mu' \in \mathcal{K}$.
	
	 Next, we show that there exists a $u$ such that $L(V_u) \leq 0$. Following \eqref{eq:optimal-u}, let 
	  \begin{equation} \label{eq:suboptimal}
	      u^\epsilon = -R^{-1}G^T\nabla_x V_u,
	  \end{equation} the control law corresponding to $V_u$. 
	  Notice that from the definition of $V_u$, $\nabla_x V_u = - \frac{\lambda}{\Psi_l} \nabla_x \Psi_l$ and $\nabla_{xx} V_u =
\frac{\lambda}{\Psi_l^2}(\nabla_x \Psi_l) (\nabla_x \Psi_l)^T - \frac{\lambda}{\Psi_l}
\nabla_{xx}\Psi_l$. So, $u^\epsilon = \frac{\lambda}{\Psi_l} R^{-1}G^T \nabla_x \Psi_l$.

Then, from \eqref{eq:LV},
  \begin{multline*}
     L(V_u) = -\frac{\lambda}{\Psi_l}(\nabla_x \Psi_l)^T
            (f + \frac{\lambda}{\Psi_l} G R^{-1}G^T \nabla_x \Psi_l)\\
          + \frac{1}{2} Tr\left(\left( \frac{\lambda}{\Psi_l^2}(\nabla_x \Psi_l )
            (\nabla_x \Psi_l)^T - \frac{\lambda}{\Psi_l} \nabla_{xx}\Psi_l \right)B \Sigma_\epsilon B\right)
  \end{multline*}
  where $\partial_t V_u = 0$ because $V_u$ is not a function of time. 
Applying the assumption in \eqref{eq:noise-assumption} and simplifying yields
  \begin{equation*}
      L(V_u)  = -\frac{\lambda}{\Psi_l}(\nabla_x \Psi_l)^T f 
	    -\frac{\lambda}{2 \Psi_l^2} (\nabla_{x}\Psi_l)^T \Sigma_t \nabla_{x}\Psi_l
                   -\frac{\lambda}{2\Psi_l} Tr\left(\left(\nabla_{xx}\Psi_l \right) \Sigma_t \right).
  \end{equation*}
From the first constraint in \eqref{eq:hjbjoin-sos},
  \begin{gather*}
     \frac{1}{\lambda} q \Psi_l - f^{T}(\nabla_{x}\Psi_l)
           -\frac{1}{2}Tr\left(\left(\nabla_{xx}\Psi_l\right)\Sigma_{t}\right) \leq 0 \implies\\
    - \frac{\lambda}{\Psi_l}(\nabla_{x}\Psi_l)^T f \leq -  
       q  + \frac{\lambda}{2 \Psi_l}Tr\left(\left(\nabla_{xx}\Psi_l\right)\Sigma_{t}\right). 
  \end{gather*}
Substituting this inequality into $L(V_u)$ and simplifying yields
  \begin{align}
	L(V_u) \leq -q -\frac{\lambda}{2 \Psi_l^2} (\nabla_{x}\Psi_l)^T \Sigma_t \nabla_{x}\Psi_l \leq 0 \label{eq:lvle0}
  \end{align}
because $q\ge 0$, $\lambda > 0$ and $\Sigma_t$ is positive semidefinite by definition. 
Since $V_u$ satisfies \cref{def:stochastic_clf}, $V_u$ is a SCLF.
\end{proof}

%The last two constraints in \eqref{eq:hjbjoin-sos} ensures that $\mu(|x|)\le V(x) \le \mu'(|x|)$. 
%Other constraints directly enforcing the condition $\mu(|x|)\le V(x) \le \mu'(|x|)$ of Definition \cref{def:stochastic_clf} are equally valid constraints that make the approximated solution $V_u$ a SCLF. Therefore, the choice of a specific type of constraints relies on the applications, and the ease computation.

\begin{corollary} \label{col:local-stable}
The suboptimal controller $u^\epsilon = -R^{-1}G^T\nabla_x V_u$ is stabilizing in probability within the domain $\Omega$.
\end{corollary}

\begin{proof}
This corollary is a direct consequence of the constructive proof of \cref{thm:lower-sclf,thm:stochastic_clf}.
\end{proof}

\begin{corollary} \label{col:local-asym-stable}
If $\Sigma_t$ is a positive definite matrix, the suboptimal controller $u^\epsilon = -R^{-1}G^T\nabla_x V_u\ $ is asymptotically stabilizing in probability within the domain $\Omega$.
\end{corollary}

\begin{proof}
This corollary is a direct consequence of the constructive proof of \cref{thm:lower-sclf,thm:stochastic_clf_asymp}.
In \eqref{eq:lvle0}, $L(V_u) < 0$ for $x \in \Omega \backslash \{0\}$ if $\Sigma_t$ is positive definite. Recall that $q$ is positive definite in the problem formulation.
\end{proof}

% In the context of the results in Corollary \cref{col:cost-upper}, the hierarchy of SOS programs yields stabilizing controllers with non-increasing suboptimality bound as the degree of polynomial approximation increases. 

%%%%%%%%%%%%%%%%%%%%%%%%%%%%%%%%%%%%%%%%%%%%%%%%%%%%%%%%%%%%%%%%%%%%%%
\subsection{Bound on the Total Trajectory Cost}

We conclude this section by showing that the expected total trajectory cost incurred by the system while 
operating under the suboptimal controller of \eqref{eq:suboptimal} can be bounded as
follows.

%\ypl{How should we think about the derivatives if the exact solution is not smooth?}

\begin{theorem}\label{thm:total-traj-cost}
Given the control law $u^\epsilon=-R^{-1}G^T \nabla_x V_u$,
  \begin{equation}
    J_u \le V_u\leq V^* -\lambda \log\left(1-\min\left\{1,\frac{\epsilon}{\eta}\right\}\right)
  \end{equation}
where $J_u = \mathbb{E}_{\omega_t} [\phi_{T}(x_{T})+\int_{0}^{T} q(x_t)+\frac{1}{2}u_{t}^{T}Ru_{t}~ dt]$, the expected cost of
the system when using the control law, $u^\epsilon$.
\end{theorem}
\begin{proof}
By It\^{o}'s formula,
	\begin{gather*}
		dV_u(x_t) = L(V_u)(x_t) dt + \nabla_x V_u(x_t) B(x_t) d\omega_t. 	
	\end{gather*}
where $L(V)$ is defined in \eqref{eq:LV}.
Then,
	\begin{equation}
		V_u(x_t)  = V_u(x_0,0) + \int^t_0 L(V_u)(x_s) ds+ \int^t_0 \nabla_xV_u(x_s) B(x_s) d\omega_s. \label{eq:temp}
	\end{equation}
Given that $V_u$ is derived from polynomial function $\Psi_l$, the integrals are well defined, and we can take the expectation of \eqref{eq:temp} to get
	\begin{align*}
		\mathbb{E}[V_u(x_t)] &= V_u(x_0,0) +\mathbb{E}\left[\int^t_0 L(V_u)(x_s) ds \right]
	\end{align*}
whereby the last term of (\cref{eq:temp}) drops out because the noise is assumed to have zero mean. The expectations of the other terms return the same terms because they are deterministic.
From \eqref{eq:lvle0}, 
\begin{align*}
	L(V_u) &\le  -q -\frac{\lambda}{2 \Psi_l^2} (\nabla_{x}\Psi_l)^T \Sigma_t \nabla_{x}\Psi_l \\
	&=-q -\frac{1}{2} \left(\nabla_{x}V_u\right)^{T}GR^{-1}G^{T}\left(\nabla_{x}V_u\right)\\
	&= -q -\frac{1}{2} ({u^\epsilon})^T R u^\epsilon
\end{align*}
where the first equality is given by the logarithmic transformation and the second equality is given by the control law $u^\epsilon=-R^{-1}G^T \nabla_x V_u$. 
% Hence,
% 	\begin{align*}
% 		&\int^t_0 L(V_u)(x_s) ds \\
% 		&\quad\le - \int^t_0  q(x_s) +\frac{1}{2} u_s^T R u_s ds.
% 	\end{align*}
	Therefore, 
	\begin{align*}
		\mathbb{E}_{\omega_t}[V_u(x_T)] &= V_u(x_0) +\mathbb{E}_{\omega_t}\left[\int^T_0 L(V_u)(x_s) ds\right] \\
		&\le V_u(x_0) -  \mathbb{E}_{\omega_t}\left[\int^T_0  q(x_s) +\frac{1}{2} (u^\epsilon_s)^T R u^\epsilon_s ds\right]  \\
		&=  V_u(x_0) - J(x_0,u^\epsilon)+\mathbb{E}_{\omega_t}[\phi(x_T)]
	\end{align*}
	where the last equality is given by \eqref{eq:cost-functional}.
Therefore, \[V_u(x_0) - J(x_0,u^\epsilon) \ge \mathbb{E}_{\omega_t}[V_u(x_T)-\phi(x_T)].\] By definition, $V_u(x_T) \geq
\phi(x_T)$ for all $x_T \in \Omega$. Thus, $\mathbb{E}_{\omega_t}[V_u(x_T)-\phi(x_T)] \geq 0$. Consequently, $V_u(x_0) - J(x_0,u^\epsilon) \geq 0$, and $V_u(x_0) \geq J(x_0,u^\epsilon)$. \Cref{thm:cost-upper} gives the second inequality in the theorem.
\end{proof}

%%% ------------------------------------------------------------------------------------------- %%%
% \section{Extensions}\label{sec:extensions}
%
% This section briefly summarizes an extension of the basic framework to a few related problems.

%%%%%%%%%%%%%%%%%%%%%%%%%%%%%%%%%%%%%%%%%%%%%%%%%%%%%%%%%%%%%%%%%%%%%%
\section{Linearly Solvable Approximations}\label{sec:extensions}

The approach presented in this paper would appear up to this point to be limited to systems that are linearly solvable, i.e., those that satisfy condition \eqref{eq:noise-assumption}.  However, the proposed methods may be extended to a system which does not satisfy these conditions by approximating the system with one that is linearly solvable.  One example is to introduce stochastic forcing into an otherwise deterministic system. 

We first construct a comparison theorem between HJB solutions to systems that share the same general dynamics, but with differing noise covariance. This comparison allows for the approximated value function of one system to bound the value function for another, providing pointwise bounds, and indeed SCLFs, for those that do not satisfy \eqref{eq:noise-assumption}.

\begin{proposition}\label{prop:linear_solve_approx}
Suppose $V^{a^*}$ is the solution to the HJB equation \eqref{eq:hjb-pde-value} with noise covariances $\Sigma_a$, and $V^b$ is a supersolution to \eqref{eq:hjb-pde-value} with identical parameters except the noise covariance $\Sigma^b$ where $\Sigma_b - \Sigma_a \succeq 0$, then $V^b \ge V^{a^*}$ for all $x \in \Omega$. 
\end{proposition}

\begin{proof}
From \cite[Def. 2.2]{crandall1992user}, $V$ is a viscosity supersolution to the HJB equation \eqref{eq:hjb-pde-value} with noise covariance $\Sigma$ if it satisfies
  \begin{equation}\label{eq:hjb_viscosity_super}
  0 \le-q - \left(\nabla_{x}V\right)^{T}f + \frac{1}{2}
           \left(\nabla_{x}V\right)^{T}GR^{-1}G^{T}\left(\nabla_{x}V\right)
                      - \frac{1}{2}Tr\left(\left(\nabla_{xx}V\right)B\Sigma B^{T}\right).
  \end{equation}
Since  $\Sigma_b - \Sigma_a \succeq 0$ the following trace inequality holds,
 \begin{equation*}
   Tr\left(\left(\nabla_{xx}V^a\right)B\Sigma_{b}B^{T}\right) 
           \ge Tr\left(\left(\nabla_{xx}V^a\right)B\Sigma_{a}B^{T}\right)\ .
  \end{equation*}
 Therefore, we have the inequality
  \begin{eqnarray*}
    0&\le & -q - \left(\nabla_{x}V^b\right)^{T}f + \frac{1}{2}
           \left(\nabla_{x}V^b\right)^{T}GR^{-1}G^{T}\left(\nabla_{x}V^b\right)- \frac{1}{2}Tr\left(\left(\nabla_{xx}V^b\right)B\Sigma_{b}B^{T}\right) \\
    &\le& -q - \left(\nabla_{x}V^b\right)^{T}f + \frac{1}{2}
           \left(\nabla_{x}V^b\right)^{T}GR^{-1}G^{T}\left(\nabla_{x}V^b\right)- \frac{1}{2}Tr\left(\left(\nabla_{xx}V^b\right)B\Sigma_{a}B^{T}\right)
  \end{eqnarray*}
  which implies that $V^b$ is in fact a viscosity supersolution to the system with noise covariance $\Sigma^a$ (i.e., $V^b$ satisfies \eqref{eq:hjb_viscosity_super} for $\Sigma^a$). As $V^b$ is a supersolution to the system with parameter $\Sigma^a$, then $V^b \ge V^{a^*}$.
\end{proof}

A particular class of such approximations arises from a deterministic HJB 
solution, which is not linearly solvable, but is approximated by one that is linearly solvable. Consider a deterministic system of the form
\begin{equation} d x_t = \left(f(x_t)+G(x_t)u_t\right)dt \label{eq:deterministic_dynamics}\end{equation} 
	with cost function
  	\begin{gather} 
         J(x,u)=\phi(x_{T})+\int_{0}^{T}q(x_t)+\frac{1}{2}u_t R u_t ~dt \label{eq:det-cost}
    \end{gather}
	where $\phi, q, R, f, G$, and the state and input domains are defined as in the stochastic problem in \cref{sec:hjb}. Then, the HJB equation is given by 
\begin{equation} \label{eq:det-hjb}
     0 =  q+\left(\nabla_{x}V\right)^{T}f-\frac{1}{2}
           \left(\nabla_{x}V\right)^{T}GR^{-1}G^{T}\left(\nabla_{x}V\right)
\end{equation}
and the optimal control is given by $u^* = -R^{-1}G^T \nabla_x V$. In general, \eqref{eq:det-hjb} is not a linear PDE. 

\begin{corollary} \label{col:upper-bound-det-v}
Let $V^*$ be the value function that solves \eqref{eq:det-hjb}, and $V^u$ be the upper bound solution obtained from \eqref{eq:hjbjoin-sos} where all parameters are the same as \eqref{eq:det-hjb} and $\Sigma_t$ is not zero. Then, $V^u$ is an upper bound for $V^*$ over the domain (i.e., $V^* \le V^u$).
\end{corollary}
\begin{proof}
A simple application of \cref{prop:linear_solve_approx}, where $\Sigma_a$ takes 
the form of a zero matrix, gives $V^* \le V^u$.
\end{proof}

Interestingly, using the solution from \eqref{eq:hjbjoin-sos} and the transformation $V_u = -\lambda \log \Psi_l$, the suboptimal controller $u^\epsilon = -R^{-1}G^T \nabla_x V_u$ is a stabilizing controller for the deterministic system \eqref{eq:deterministic_dynamics} if a simple condition is satisfied. This fact is shown using the Lyapunov theorem for deterministic systems introduced next \cite{sontag1983lyapunov}.

\begin{definition} \label{def:clf}
	Given the system \eqref{eq:deterministic_dynamics} and cost function \eqref{eq:det-cost}, a control Lyapunov function (CLF) is a proper positive definite function $\mathcal{V} \in \mathcal{C}^1$ on a compact domain $\Omega \cup \{0\} $ such that 
	\begin{gather}
		\mathcal{V}(0) = 0,~ \mathcal{V}(x) \ge \mu(|x|)\quad \forall~ x \in \Omega\backslash \{0\} \label{eq:clf}\\
		\exists~ u(x) ~ s.t. ~ (\nabla_x \mathcal{V})^T (f+Gu) \le 0 \quad \forall~ x \in \Omega\backslash \{0\}  \nonumber
	\end{gather}
	where $\mu \in \mathcal{K}$.
\end{definition}
\begin{theorem}[\cite{sontag1983lyapunov} Thm. 2.5]\label{thm:clf-lyapunov}
	Given a system \eqref{eq:deterministic_dynamics} and cost function \eqref{eq:det-cost}, if there exists a CLF $V$ and a $u$ satisfying \cref{def:clf}, then the controlled system is stable, and $u$ is a stabilizing controller. Furthermore, if $(\nabla_x V)^T (f+Gu) < 0$ for all  $x \in \Omega\backslash \{0\}$, the controlled system is asymptotically stable, and $u$ is an asymptotically stabilizing controller.
\end{theorem}

Verifying that the controller $u^\epsilon = -R^{-1}G^T \nabla_x V_u$ is in fact stabilizing and that $V_u$ is a CLF may be seen as follows. 

\begin{corollary} \label{col:det-lyap}
Given the controller $u^\epsilon = -R^{-1}G^T \nabla_x V_u$, if 
\begin{equation*}Tr \left( \left( \nabla_{xx} V_u \right)\Sigma_t \right) \ge 0 \quad \forall ~x \in \Omega \backslash \{0\}, \end{equation*}
then $u^\epsilon$ is a stabilizing controller for \eqref{eq:deterministic_dynamics}. If 
\begin{equation*}Tr \left( \left( \nabla_{xx} V_u \right) \Sigma_t \right) > 0 \quad \forall~ x \in \Omega \backslash \{0\}, \end{equation*}
then $u^\epsilon$ is an asymptotically stabilizing controller for \eqref{eq:deterministic_dynamics}.
\end{corollary}

\begin{proof}
Recall that from the proof of \cref{thm:lower-sclf}, all conditions in \cref{def:clf} are satisfied by $V_u$ except \eqref{eq:clf}. 
To show that $V_u$ satisfies \eqref{eq:clf}, rearrange \eqref{eq:hjb-pde-value} to yield the following
  \begin{align*}
      (\nabla_x V_u )^T (f+Gu^\epsilon)&=(\nabla_x V_u)^T f -(\nabla_x V_u)^T G R^{-1} G^T (\nabla_x V_u) \\
 & \le  -q -\frac{1}{2} (\nabla_x V_u)^T G R^{-1} G^T (\nabla_x V_u) - \frac{1}{2} Tr \left( \left( \nabla_{xx} V_u \right)\Sigma_t \right)
  \end{align*}
  where $\Sigma_t = B \Sigma_\epsilon B^T$.
Recall that $q$ and $R$ are positive definite. If $Tr \left( \left( \nabla_{xx} V_u \right)\Sigma_t \right) \ge 0$ for all $x \in \Omega\backslash \{0\}$, then $(\nabla_x V_u)^T (f+Gu^\epsilon) \le 0$ implying that $V^u$ is a CLF and $u^\epsilon$ is a stabilizing controller by \cref{thm:clf-lyapunov}. Furthermore,  if $Tr \left( \left( \nabla_{xx} V_u \right) \Sigma_t \right) > 0$ for all $x \in \Omega\backslash \{0\}$, $u^\epsilon$ is an asymptotically stabilizing controller.
\end{proof}

In the deterministic case, $\Sigma_t$ is free variable that can be chosen to be small according to the equality \eqref{eq:noise-assumption}. Hence, \eqref{eq:noise-assumption} is no longer a constraint or an assumption, but it serves as a design principle for obtaining a CLF for system \eqref{eq:deterministic_dynamics}. 
Furthermore, given a $\Sigma_t$, the trace condition in \cref{col:det-lyap} is easily enforced in \eqref{eq:hjbjoin-sos} by adding one extra constraint in the optimization problem. Thus, the optimization problem \eqref{eq:hjbjoin-sos} can also produce a CLF for the corresponding deterministic system, with analytical results from the \cref{sec:analysis}, including a priori trajectory suboptimality bounds (\cref{thm:total-traj-cost}), inherited as well.

\setcounter{figure}{0}
 \begin{figure}[b!]%[b!]%[h!]
   \centering
   \includegraphics[trim = 12mm 65mm 15mm 70mm, clip, width=0.45\textwidth]{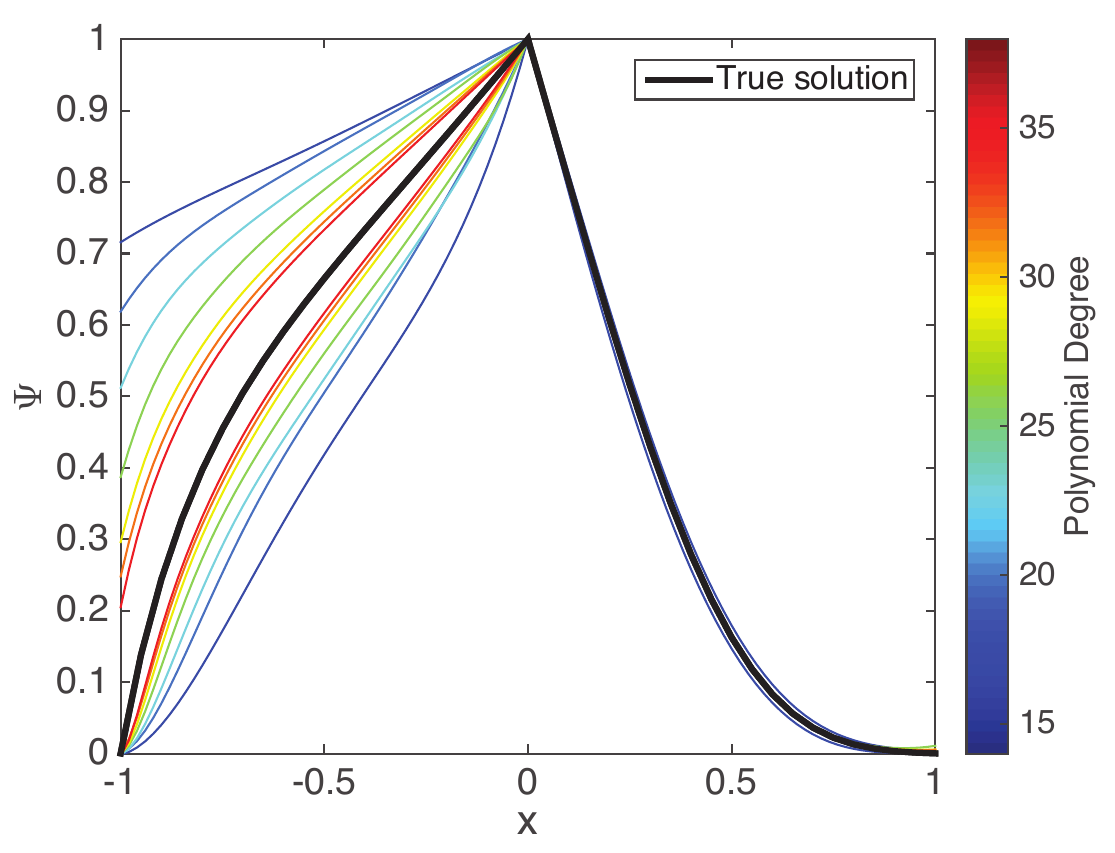}
   \caption{The desirability function of system \eqref{eq:1dexample} for varying polynomial degree. The true solution is the black curve.}
      \label{fig:Desirability-ex-1}
 \end{figure}

%%% ------------------------------------------------------------------------------------------- %%%
\section{Numeric Examples}\label{sec:simulation}

This section studies the computational characteristics of this method using two examples -- a scalar system and a two-dimensional system. In the following problems, the optimization parser YALMIP \cite{yalmip} was used in conjunction with the semidefinite optimization package MOSEK
\cite{andersen2000mosek}. In both examples, the continuous system is integrated numerically using Euler integration with step size of 0.005s during simulations.

 % %%%%%%%%%%%%%%%%%%%%%%%%%%%%%%%%%%%%%%%%%%%%%%%%%%%%%%%%%%%%%%%%%%%%%%
%  \subsection{Scalar Stable System Example}
%  Consider the following stable nonlinear system
%    \begin{equation*} %\label{eq:ex_1}
%        dx=\left(x^{3}+x^{2}-4x+u\right)dt+d\omega
%    \end{equation*}
%  on the domain $x\in \Omega = \{x \mid -1 \le x \le 1 \}$. The goal is to stabilize the system at the origin such that
%  $\phi(-1)=20e^{-10}$, $\phi(1)=20e^{-10}$ and $\phi(0) = 1$. We set $\Sigma_\epsilon=1$, $G=1$, $B=1$,
%  $q=10$, and $R=1$. In this example, instead of the trace constraint in \eqref{eq:hjbjoin-sos}, we use the constraints $\frac{d \Psi}{dx} \ge 0 ~\forall x < 0$ and $\frac{d \Psi}{dx} \le 0 ~\forall x \ge 0$ to simplify computation. In addition, because this example is a scalar system, we can divide the domain $\Omega$ into two partitions $x\le 0$ and $x \ge 0$ and solve them simultaneously.
%
%  The desirability function from solving \eqref{eq:hjbjoin-sos} for varying
%  polynomial degrees are shown in Figure \cref{fig:Desirability-ex-1}. The kink at the origin is expected because the HJB PDE solution is always smooth within a domain, but it is not necessarily smooth at the boundary. In this example, the origin is a boundary. The approximation error $\epsilon$ for both partitions is shown in Figure \cref{fig:ex-1-solution-epsilon} for increasing polynomial degree. As seen in the plots, the approximation improves as the degree of polynomial increases.
%

 %%%%%%%%%%%%%%%%%%%%%%%%%%%%%%%%%%%%%%%%%%%%%%%%%%%%%%%%%%%%%%%%%%%%%%

 \subsection{Scalar Unstable System} 
Consider the following scalar unstable nonlinear system
   \begin{equation} 
       dx=\left(x^{3}+5x^{2}+x+u\right)dt+d\omega \label{eq:1dexample}
   \end{equation}
 on the domain $x\in \Omega = \{x \mid -1 \le x \le 1 \}$. The noise model considered is Gaussian white noise with zero mean and variance $\Sigma_\epsilon=1$. The goal is to stabilize the system at the origin. We choose the boundary at two ends of the domain to be $\Psi(-1)=20e^{-10}$ and $\Psi(1)=20e^{-10}$. At the origin, the boundary is set as $\Psi(0) = 1$. We set $q=x^2$, and $R=1$. 
In the one dimensional case, the origin, which is a boundary, divides the domain into two partitions, $x\le 0$ and $x \ge 0$. 
Because of the natural division of the domain, the solutions for both domains can be represented by smooth polynomial respectively, and solved independently. The simulation is terminated when the trajectories enter the interval $[-0.005,0.005]$ centered on the origin. 

The desirability functions that result from solving \eqref{eq:hjbjoin-sos} for varying polynomial degrees are shown in \cref{fig:Desirability-ex-1}. The true solution is computed by solving the HJB directly in Mathematica \cite{mathematica}. 
The kink at the origin is expected because the HJB PDE solution is not necessarily smooth at the boundary, and in this instance the origin is a zero-cost boundary. 

The approximation error $\epsilon$ for both partitions is shown in \cref{fig:ex-1-solution-epsilon} for increasing polynomial degree. As seen in the plots, the approximation improves as the polynomial degree increases. Polynomial degrees below 14 are not feasible, hence this data is absent in the plots. The suboptimal solution converges faster for $x>0$ than for $x<0$ when the degree of polynomial increases because the true solution for $x>0$ has a simple quadratic-like shape that can be easily represented as a low degree SOS function. 

\Cref{fig:ex-1-solution-traj} shows sample trajectories using the controller computed from optimization problem \eqref{eq:hjbjoin-sos} for different polynomial degrees. The controllers are stabilizing for six randomly chosen initial points. Unsurprisingly, the suboptimal solutions with low pointwise error result in the system converging towards the origin faster.

To compare between $J_u$ and $V_u$, a Monte Carlo experiment is illustrated in \cref{fig:ex-1-solution-cost}. For each polynomial degree that is feasible, the controller obtained from $\Psi_l$ in optimization problem \eqref{eq:hjbjoin-sos} is implemented in 30 simulations of the system subject to random samples of Gaussian white noise with $\Sigma_\epsilon = 1$. The initial condition is fixed at $x_0 =-0.5$. In the figure, $V^u \geq J^u$ as expected, and the difference between the two decreases with increasing $d$.

% Figure \cref{fig:ex-1-solution-cost} shows the comparison between $J_u$ and $V_u$ for different polynomial degrees whereby $J_u$ is the expected cost and $V_u$ is the value function computed from $\Psi_l$ in optimization problem \eqref{eq:hjbjoin-sos}. 

 \setcounter{figure}{1}
 \begin{figure*}[t!]
       \centering
 	  \subfigure[ ]{\includegraphics[trim =12mm 65mm 15mm 70mm, clip,width=0.32\textwidth]{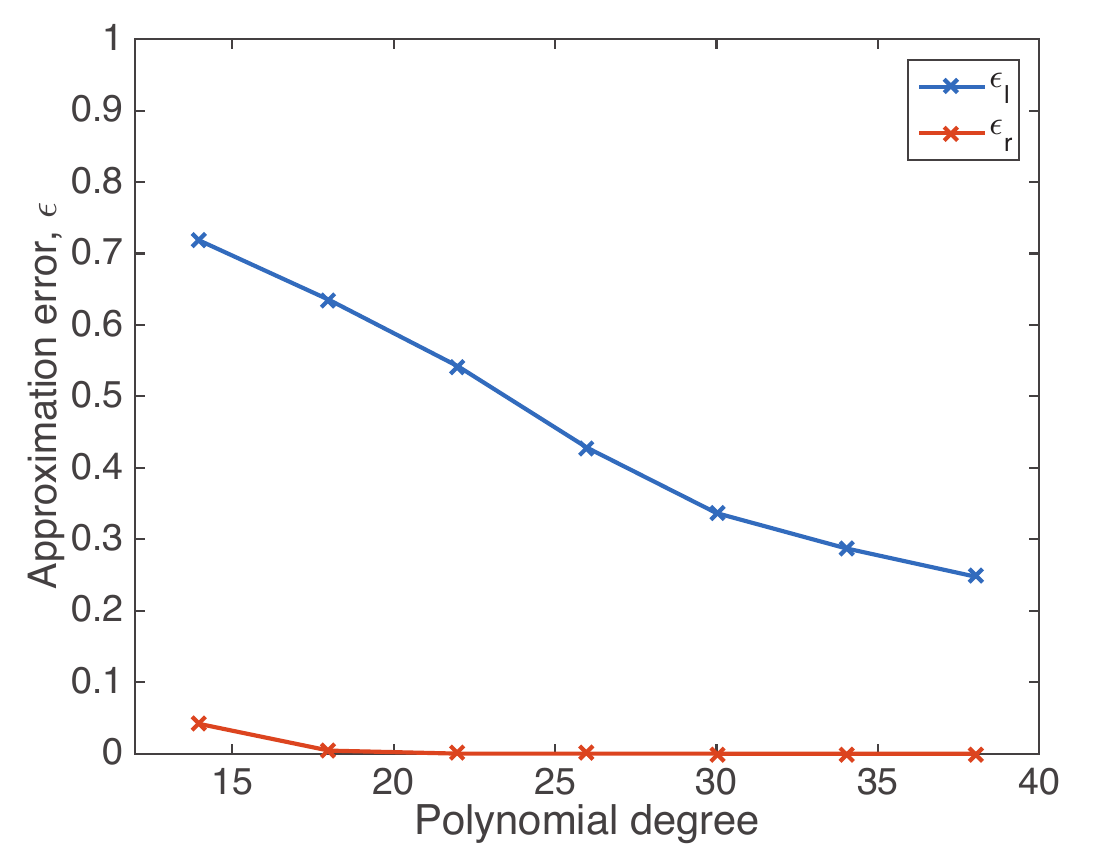}\label{fig:ex-1-solution-epsilon}  }  
 	  \subfigure[ ]{\includegraphics[trim =12mm 65mm 15mm 70mm, clip,width=0.32\textwidth]{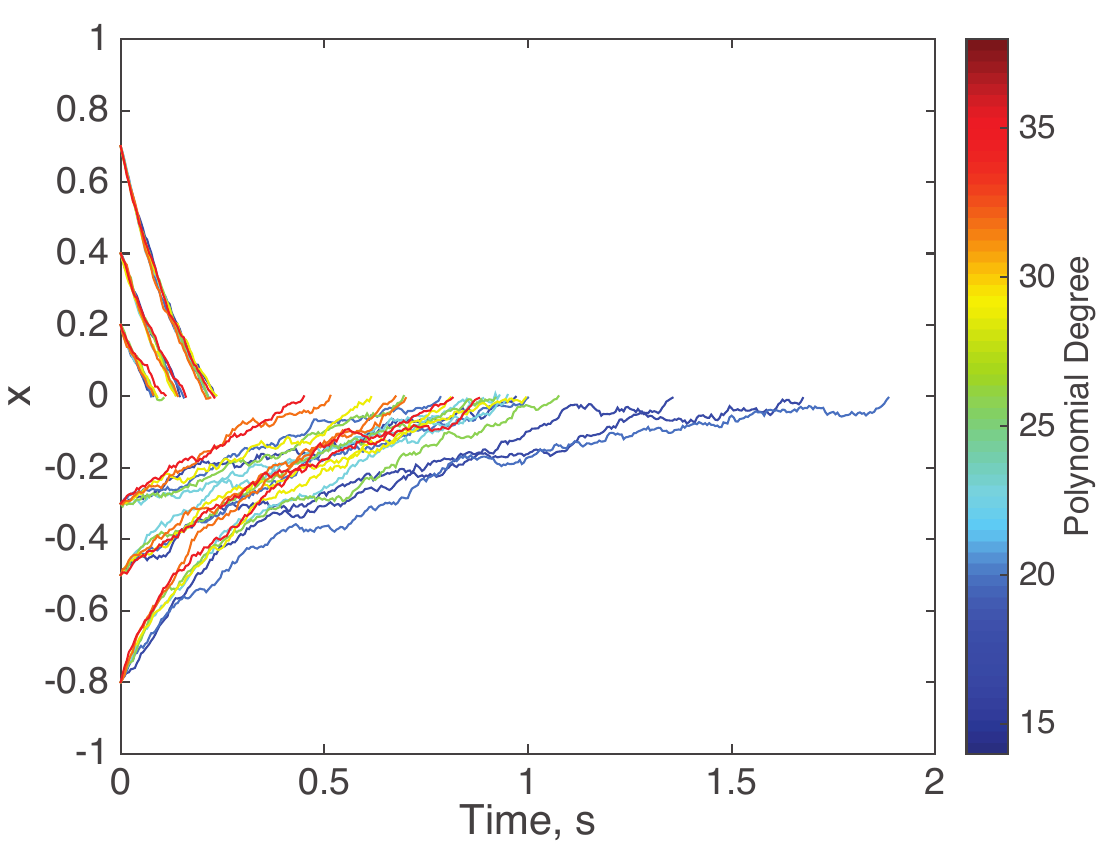} \label{fig:ex-1-solution-traj}}  
 	  \subfigure[ ]{\includegraphics[trim =15mm 65mm 15mm 70mm, clip,width=0.32\textwidth]{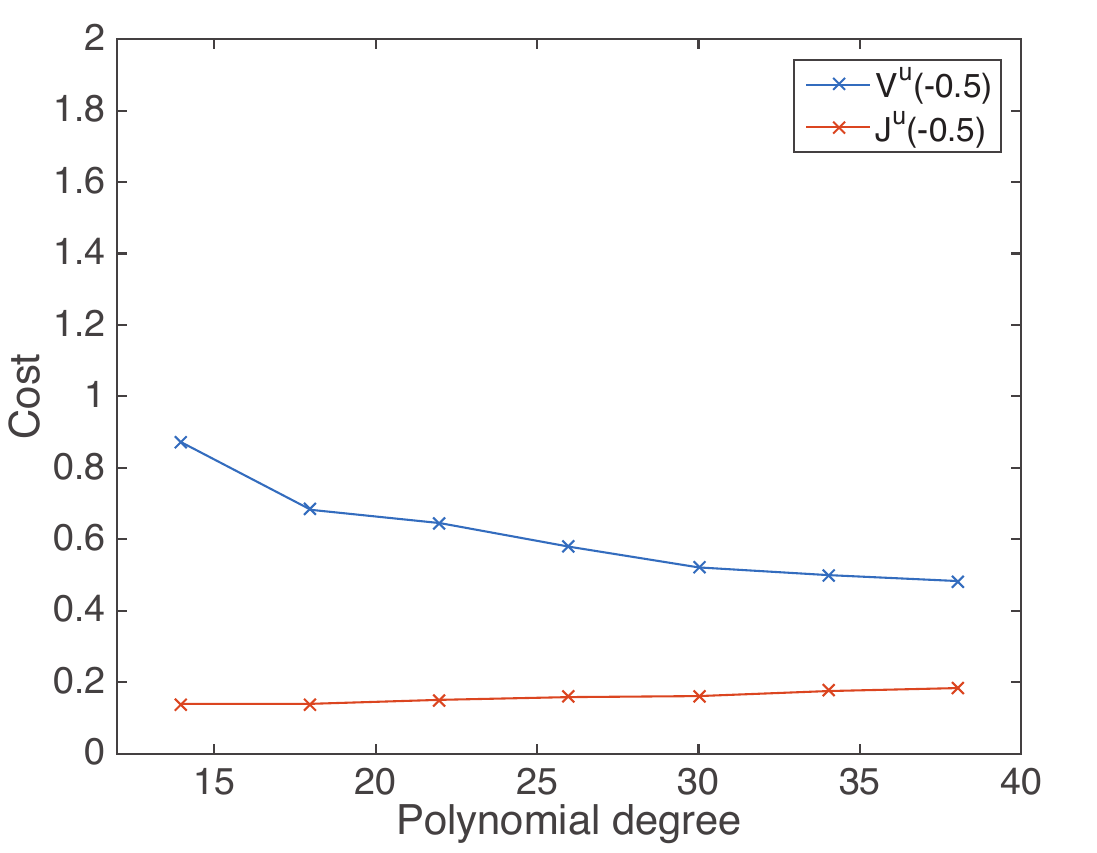}\label{fig:ex-1-solution-cost}} 	  
       \vspace{-.2cm}
       \caption{Computational results of system \eqref{eq:1dexample}. (a) Convergence of the objective function of  \eqref{eq:hjbjoin-sos} as the degree of polynomial increases. The approximation error for $x \le 0$ is denoted as $\epsilon_l$ and the approximation error for $x \ge 0$ is denoted as $\epsilon_r$. (b) Sample trajectories using controller computed from optimization problem \eqref{eq:hjbjoin-sos} with different polynomial degrees starting from six randomly chosen initial points. (c) The comparison between $J_u$ and $V_u$ for different polynomial degrees whereby $J_u$ is the expected cost and $V_u$ is the value function computed from optimization problem \eqref{eq:hjbjoin-sos}. The initial condition is fixed at $x_0 = -0.5$.}     
  \end{figure*}

%%%%%%%%%%%%%%%%%%%%%%%%%%%%%%%%%%%%%%%%%%%%%%%%%%%%%%%%%%%%%%%%%%%%%
 \subsection{Two Dimensional System}

In the following example, we demonstrate the power of this technique on a 2-dimensional system. Consider a nonlinear 2-dimensional problem example with the following dynamics:
  \begin{equation}
   \left[\begin{array}{c}
    dx\\ 
    dy
   \end{array}\right] =  \left(2\left[\begin{array}{c}
                             x^5-x^{3}-x+x y^4\\
                            y^5-y^{3}-y+y x^4
   \end{array}\right]+\left[\begin{array}{c}
    x~ u_{1}\\
    y ~u_{2}
  \end{array}\right]\right)dt
       + \left[\begin{array}{c}
     x~ d\omega_{1}\\
     y ~d\omega_{2}
  \end{array}\right]. \label{eq:2dexample}
\end{equation}
The goal is to reach the origin at the boundary of the domain $\Omega = \{(x,y) \mid -1\le x \le 1, -1 \le y \le 1\}$. The control penalty is $R=I_{2\times 2}$, and state cost is $q(x)=x^2+y^2$. The boundary conditions for the sides at $x=1,x=-1,y=1$, and $y=-1$ are set to $\phi(x,y)=5$, while at the origin, the boundary has cost $\phi(0,0)= 0$. The noise model considered is Gaussian white noise with zero mean and an identity covariance matrix. 

 \setcounter{figure}{2}
 \begin{figure}[t!]%[b!]
 \centering
 \subfigure[$\Psi$, Degree = 10]{\includegraphics[trim = 12mm 65mm 15mm 70mm, clip, width=0.23\textwidth]{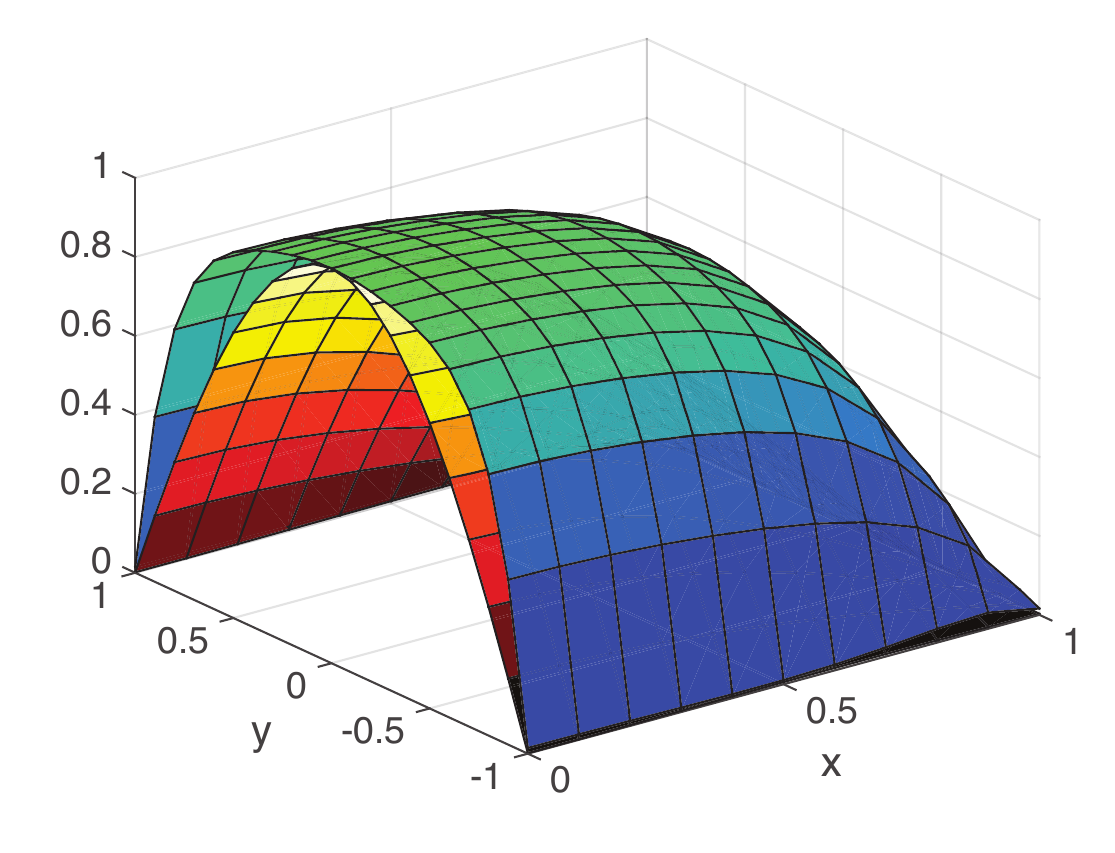}}
 \subfigure[$\Psi$, Degree = 20]{\includegraphics[trim = 12mm 65mm 15mm 70mm, clip, width=0.23\textwidth]{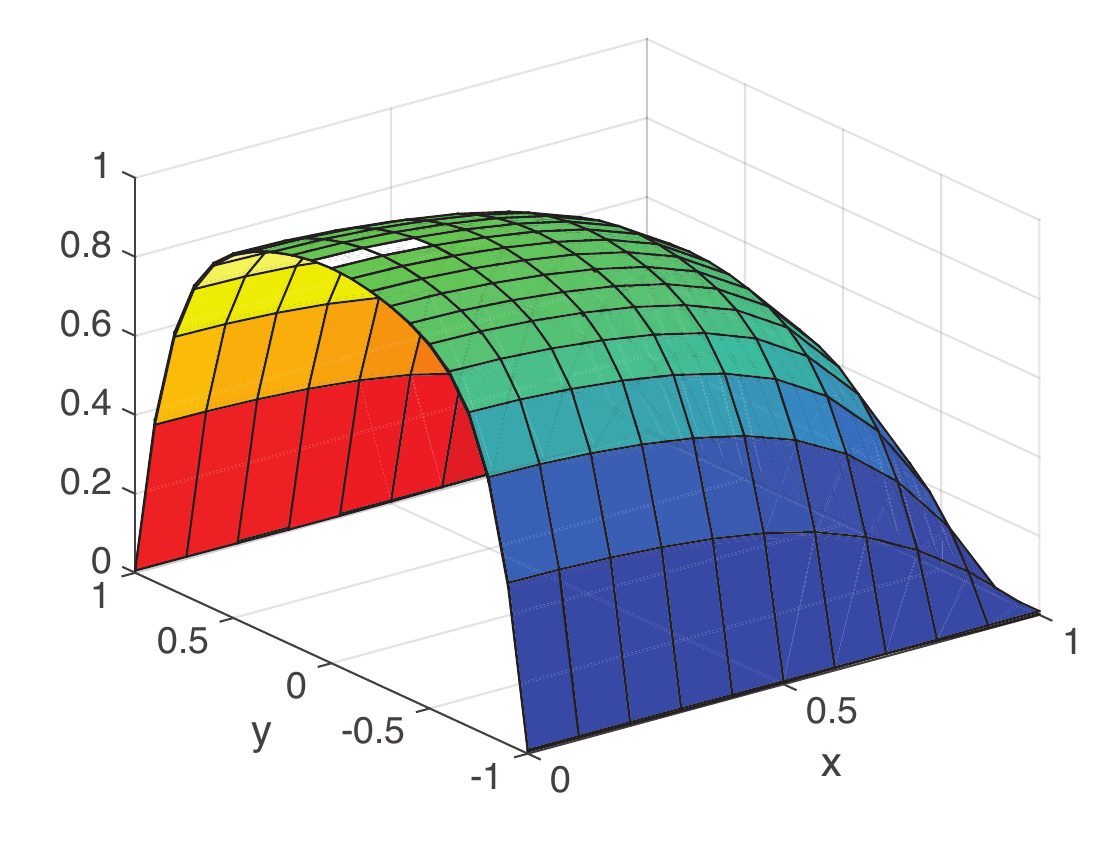}}
 \subfigure[$V$, Degree = 10]{\includegraphics[trim = 12mm 65mm 15mm 70mm, clip, width=0.23\textwidth]{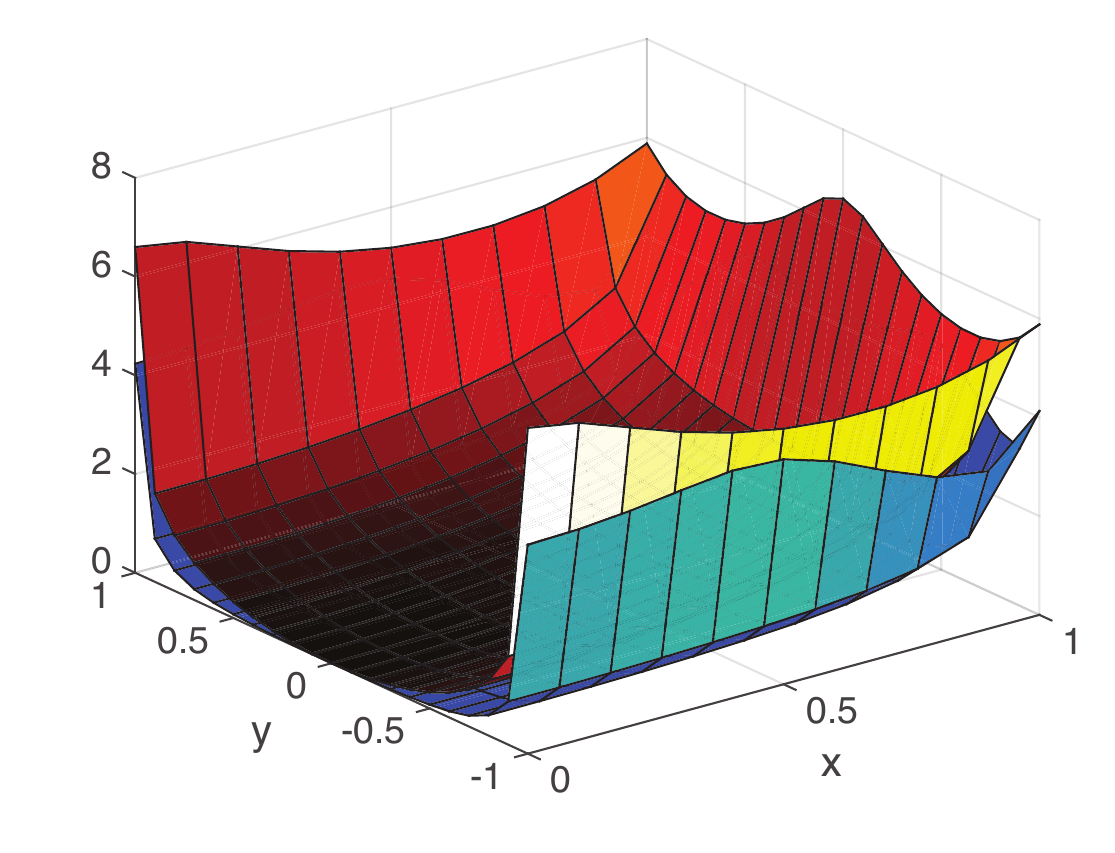}}\label{fig:Results-of-multidimensional-v10}
 \subfigure[$V$, Degree = 20]{\includegraphics[trim = 12mm 65mm 15mm 70mm, clip, width=0.23\textwidth]{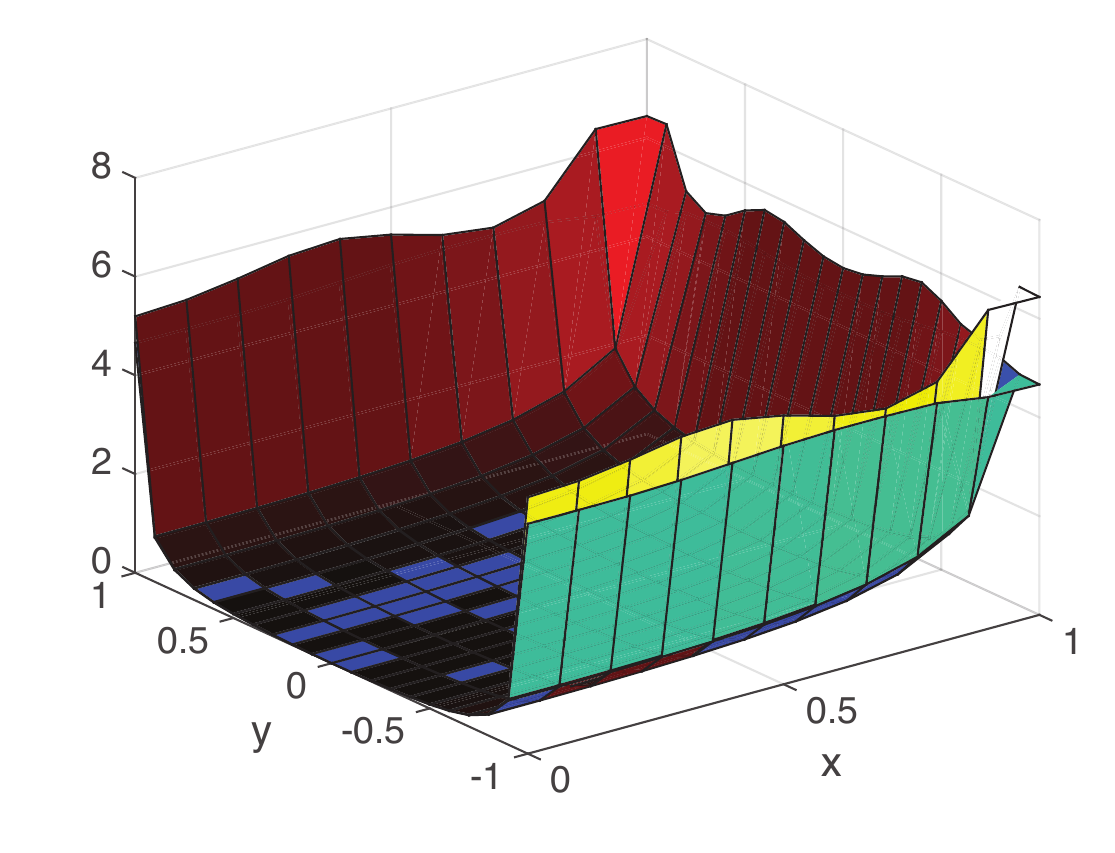}}\label{fig:Results-of-multidimensional-v20}
 \caption{Approximated desirability functions and value functions for \eqref{eq:2dexample} when polynomial degrees are 10 and 20. In (a) and (b), the blue sheets are the upper bound solutions $\Psi_u$ and the red sheets are the lower bound solutions $\Psi_l$. The corresponding value functions are shown in (c) and (d) respectively.}
 \label{fig:Results-of-multidimensional}
 \end{figure}

 \setcounter{figure}{3}
 \begin{figure*}[b!]%[t!]%[h!]
      \centering
	  \subfigure[ ]{\includegraphics[trim =12mm 65mm 15mm 70mm, clip,width=0.32\textwidth]{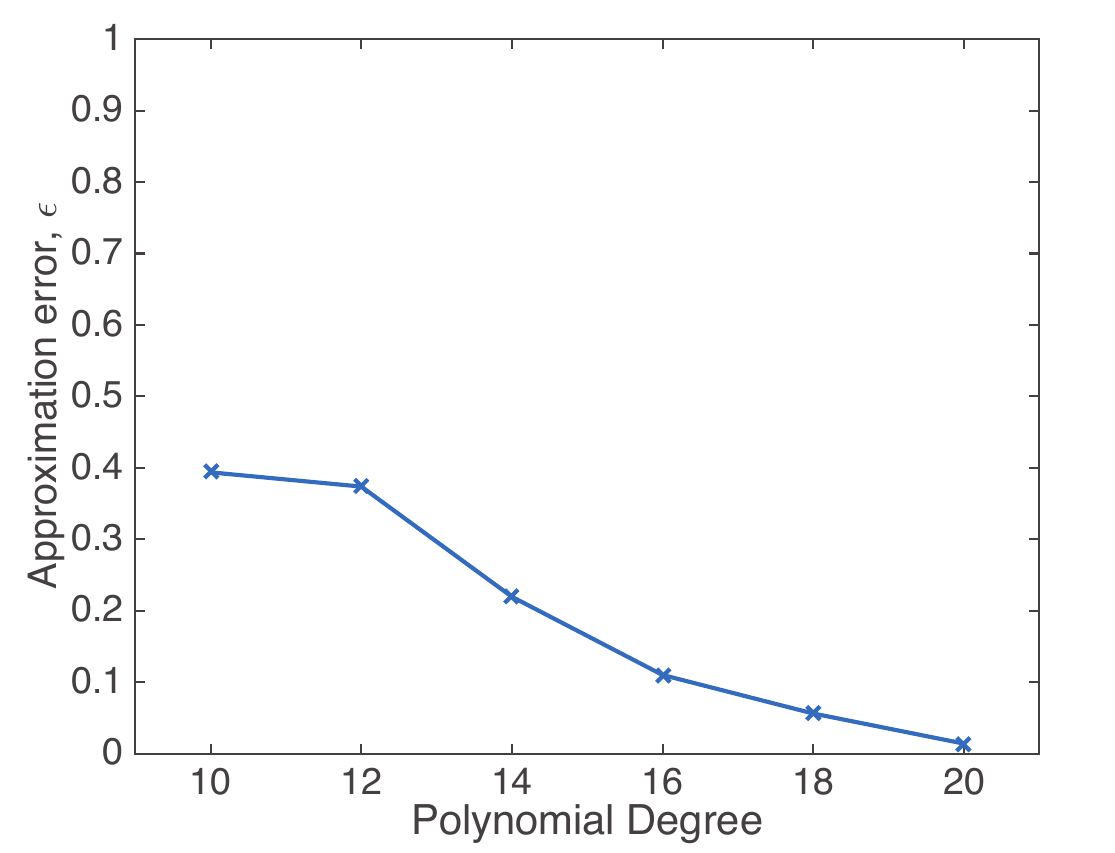}\label{fig:ex-2-solution-epsilon} }     	 
	  \subfigure[ ]{\includegraphics[trim =12mm 65mm 15mm 70mm, clip,width=0.32\textwidth]{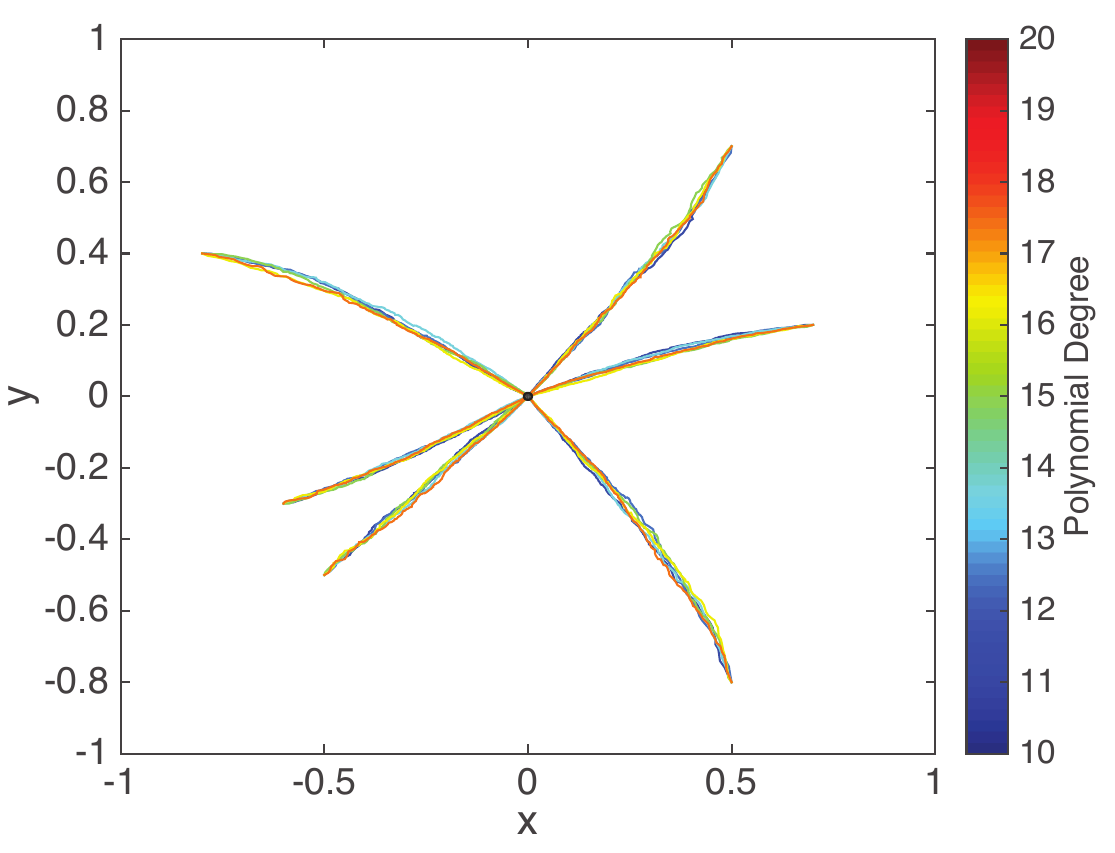}\label{fig:ex-2-solution-traj}} 
	   \subfigure[ ]{\includegraphics[trim =15mm 65mm 15mm 70mm, clip,width=0.32\textwidth]{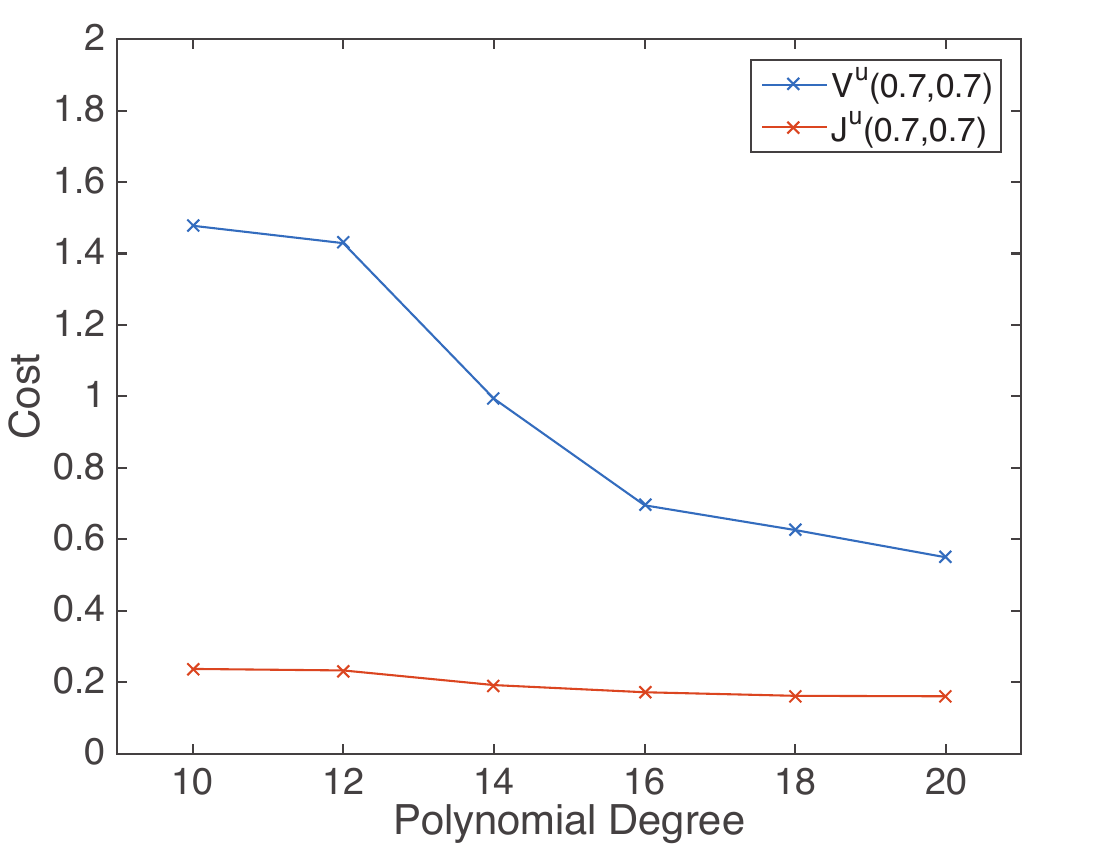}\label{fig:ex-2-solution-cost}} 
	   \vspace{-.2cm}
      \caption{Computational results of system \eqref{eq:2dexample}. (a) Convergence of the variables in the objective function of \eqref{eq:hjbjoin-sos}. (b) Sample trajectories using controller from optimization problem \eqref{eq:hjbjoin-sos} with different polynomial degrees starting from six randomly chosen initial points. (c) The comparison between $J_u$, the expected cost, and $V_u$ the value function for different polynomial degrees from optimization problem \eqref{eq:hjbjoin-sos}. The initial condition is fixed at $x_0 = (0.7, 0.7)$. }     
 \end{figure*}

The approximated desirability functions and their corresponding value functions are shown in \cref{fig:Results-of-multidimensional}, with half of the domain $x \in [0,1]$ shown in order to view the gaps between the upper and lower bound solutions. %When the polynomial degree is 20, the upper and lower bound solutions are numerically identical in many regions. 
% Hence, the blue and red sheets in Figure \cref{fig:Results-of-multidimensional-v10} and \cref{fig:Results-of-multidimensional-v20} crosses each other. 
\Cref{fig:ex-2-solution-epsilon} shows the convergence of the objective function of optimization problem \eqref{eq:hjbjoin-sos} as the degree of polynomial increases. There is no data below degree of 10 because the optimization problem is not feasible in these cases. As shown in \cref{fig:ex-2-solution-traj}, sample trajectories starting from six different initial points shows that the controllers computed from $\Psi_l$ for various degrees arrive at the origin. The trajectory is considered at the origin if it is within a distance of 0.01 from the origin.  

Similar to the scalar example, a Monte Carlo experiment is performed to compare between $J_u$ and $V_u$. For each polynomial degree that is feasible, the controller obtained from $\Psi_l$ in optimization problem \eqref{eq:hjbjoin-sos} is implemented in 30 simulations of the system subject to random samples of Gaussian white noise with $\Sigma_\epsilon = I_{2\times 2}$. The initial condition is fixed at $x_0 =(0.7,0.7)$. \Cref{fig:ex-2-solution-cost} shows the comparison between $J_u$ and $V_u$ for different polynomial degrees whereby $J_u$ is the expected cost and $V_u$ is the value function computed from $\Psi_l$ in optimization problem \eqref{eq:hjbjoin-sos}. As expected, $V^u \geq J^u$.

%%%%%%%%%%%%%%%%%%%%%%%%%%%%%%%%%%%%%%%%%%%%%%%%%%%%%%%%%%%%%%%%%%%%%%
\section{Conclusion}\label{sec:conclusion}

This paper has proposed a new method to approximate the solution to a class of optimal
control problems for stochastic nonlinear systems via SOS programming.
Analytical results provide guarantees on the suboptimality of trajectories when using the
approximate solutions for controller design. Consequently, one can synthesize a suboptimal stabilizing
controller for a large class of stochastic nonlinear dynamical systems.

As is commonly seen when using SOS programming, the numerics of the SDP may be cumbersome in practice. There are a number of avenues for future work aimed at improving the practical performance. First, the monomials of the polynomial approximation can be chosen strategically in order to decrease computation time while achieving high accuracy. A promising future direction is the synthesis of the work presented here with that of \cite{horowitz2014linear}, wherein the curse of dimensionality is avoided via the strategic choice of basis functions. To improve the numerical conditioning of these optimization techniques, a domain partitioning technique is studied in \cite{horowitz2014admm}, wherein the alternating direction method of multipliers is used to enable both parallelization and a solution representation that varies in resolution over the domain. In addition, there exists a growing body of literature towards increasing the numeric stability and scalability of SOS techniques \cite{permenter2014basis,ahmadi2014towards}. %The incorporation of these techniques into the present work is under investigation.

\bibliographystyle{siamplain} 
\bibliography{references}

% \begin{IEEEbiographynophoto}{Yoke Peng Leong}
% received her B.S. and M.S. degrees in mechanical engineering from Northwestern University in 2012. She is currently pursuing the Ph.D degree in Control and Dynamical Systems at California Institute of Technology.
% \end{IEEEbiographynophoto}
%
% \begin{IEEEbiographynophoto}{Matanya B. Horowitz}
% Biography text here.
% \end{IEEEbiographynophoto}
%
% \begin{IEEEbiographynophoto}{Joel W. Burdick}
% Biography text here.
% \end{IEEEbiographynophoto}

\end{document}